\documentclass[review,1p,times,authoryear]{elsarticle}

\usepackage{graphicx}      
\usepackage{amsmath}
\usepackage{amssymb}
\usepackage[english]{babel}
\usepackage{algorithm}
\usepackage{algpseudocode}
\usepackage{amsthm}
\usepackage{subcaption}
\usepackage{float}
\usepackage{framed}
\usepackage{tabto}
\usepackage{hyperref}

\usepackage[final]{changes}

\usepackage{xpatch}
\usepackage{nomencl} 
\makenomenclature
\xpatchcmd{\thenomenclature}{%
  \section*{\nomname}
}{
}{\typeout{Success}}{\typeout{Failure}}

\usepackage{etoolbox}
\renewcommand\nomgroup[1]{%
  \item[\bfseries
  \ifstrequal{#1}{M}{Main Symbols}{%
  \ifstrequal{#1}{S}{Subscripts and Superscripts}{%
  \ifstrequal{#1}{A}{Acronyms}{}%
  }}%
]}

\usepackage{multicol}
\makeatletter
\@ifundefined{chapter}
  {\def\wilh@nomsection{section}}
  {\def\wilh@nomsection{chapter}}

\def\thenomenclature{%
  \begin{multicols}{2}[%
    \csname\wilh@nomsection\endcsname*{\nomname}
    \if@intoc\addcontentsline{toc}{\wilh@nomsection}{\nomname}\fi
    \nompreamble]
  \list{}{%
    \labelwidth\nom@tempdim
    \leftmargin\labelwidth
    \advance\leftmargin\labelsep
    \itemsep\nomitemsep
    \let\makelabel\nomlabel}%
}
\def\endthenomenclature{%
  \endlist
  \end{multicols}
  \nompostamble}
\makeatother

\DeclareMathOperator*{\argmin}{arg\,min}

\definechangesauthor[name={Benoit Chachuat}, color=orange]{BC}
\setaddedmarkup{\uwave{#1}}
\setdeletedmarkup{\color{red}\sout{#1}}

\newtheorem{assumption}{Assumption}[section]

\journal{Computers \& Chemical Engineering}

\begin{document}

\begin{frontmatter}

\title{\replaced{Real-Time Optimization}{Modifier Adaptation} Meets Bayesian Optimization and Derivative-Free Optimization\added{: A Tale of Modifier Adaptation}}

\author[ICL]{E. A. del Rio Chanona\fnref{first_author}}
\author[UCL]{P. Petsagkourakis\fnref{first_author}}
\fntext[first_author]{Equal contributors}
\author[NTNU]{E. Bradford}
\author[USP]{J. E. Alves Graciano\fnref{new_address}}
\fntext[new_address]{Present address: Radix Engenharia e Software, Rio de Janeiro, Brazil}
\author[ICL]{B.~Chachuat\corref{corraut}}
\cortext[corraut]{Corresponding author}
\ead{b.chachuat@imperial.ac.uk}

\address[ICL]{Centre for Process Systems Engineering, Department of Chemical Engineering, Imperial College London, UK}
\address[UCL]{Centre for Process Systems Engineering, Department of Chemical Engineering, University College London, UK}
\address[NTNU]{Department of Engineering Cybernetics, Norwegian University of Science and Technology, Trondheim, Norway}
\address[USP]{Universidade de S\~ao Paulo, Escola Politecnica, Departamento de Engenharia Qu\'imica, S\~ao Paulo, Brazil}


\begin{abstract}
This paper investigates a new class of modifier-adaptation schemes to overcome plant-model mismatch in real-time optimization of uncertain processes. The main contribution lies in the integration of concepts from the \replaced{fields}{areas} of Bayesian optimization and derivative-free optimization. The proposed schemes embed a physical model and rely on trust-region ideas to minimize risk during the exploration, while employing Gaussian process regression to capture the plant-model mismatch in a non-parametric way and drive the exploration by means of acquisition functions. The benefits of using an acquisition function, knowing the process noise level, or specifying a nominal process model are \replaced{analyzed}{illustrated} on numerical case studies, including a semi-batch photobioreactor optimization problem \added{with a dozen decision variables}.
\end{abstract}

\begin{keyword}
real-time optimization \sep modifier adaptation \sep trust region \sep Gaussian process regression \sep Bayesian optimization \sep acquisition function \sep model-free RTO
\end{keyword}

\end{frontmatter}


\section{Introduction}

The business benefits of real-time optimization (RTO) in the oil-and-gas and chemical sectors are not disputed \citep{Darby2011, Camara2016}. Despite this, the deployment and penetration  of this technology have remained relatively low. The causes for this are many, but in particular, companies invariably need to employ highly-qualified process control engineers to design, install and continually maintain RTO systems to preserve benefits. These systems rely on knowledge-driven (mechanistic) models, and in those processes where the optimization execution period is much longer than the closed-loop process dynamics, steady-state models are commonly employed to conduct the optimization \citep{Marlin1997}. Traditionally, the model is updated in real-time using process measurements, before repeating the optimization on a time-scale of hours to days. This two-step RTO scheme, often referred to as model-adaptation strategy, is both intuitive and popular but it can hinder convergence to a plant's optimal operating point due to lack of integration between the model-update and optimization steps, especially in the presence of plant-model mismatch  \citep{Tatjewski2002,Gao2005,Tejeda2019}. This has fueled the development of alternative adaptation paradigms in RTO \citep{Engell2007,Chachuat2009}, such as modifier adaptation \citep{Marchetti2009}.

\nomenclature[A]{RTO}{real-time optimization}
\nomenclature[A]{MA}{modifier adaptation}

Similar to the two-step RTO scheme, modifier adaptation embeds the available process model into a nonlinear optimization problem that is solved on every RTO execution. The key difference is that the process measurements are now used to update the so-called modifiers that are added to the cost and constraint functions in the optimization model, while keeping a nominal process model. This methodology greatly alleviates the problem of offset from the actual plant optimum, by ensuring that the KKT conditions determined by the model match those of the plant upon convergence \citep{Marchetti2009}. However, this desideratum comes at the cost of having to estimate the cost and constraint gradients from process measurements.

Inferring gradient information from noisy process measurements is challenging, but nonetheless key to the effectiveness and reliability of modifier adaptation \citep{Bunin2013,Jeong2018}. Variants of the modifier-adaptation principle in order to mitigate this burden are surveyed by \citet{Marchetti2016}. They include recursive update schemes that exploit past steady-state operating points \citep{Gao2005,Marchetti2010,Rodger2011}, selective adaptation schemes that rely on directional derivatives \citep{Costello2016}, as well as schemes that take advantage of transient process measurements \citep{Francois2014,Krishnamoorthy2018,Speakman2020}. Other variants do not require estimating plant gradients explicitly. The nested modifier-adaption scheme by \citet{Navia2015} embeds the modified optimization model into an outer problem that optimizes over the gradient modifiers using a derivative-free algorithm. \citet{Gao2016} proposed to combine quadratic surrogates trained on available plant data with a nominal mechanistic model in order to account for curvature information and filter out the process noise. Likewise, \citet{Singhal2016} investigated data-driven approaches based on quadratic surrogates as modifiers for the predicted cost and constraint functions and devised an online adaptation strategy for the surrogates inspired by trust-region ideas. More recently, \citet{Ferreira2018} were the first to consider Gaussian processes (GPs), trained from past measurement information, as the cost and constraint modifiers. \deleted{Then }\citet{delRio2019} developed this strategy further by introducing modifier-adaptation schemes that rely on trust regions to capture the GPs' ability to capture the cost and constraint mismatch. \added{Lately, \citet{Shukla2020} investigated convergence certificates for such schemes and confirmed the benefits of using GP surrogates owing to their probabilistic full-linearity properties.} But the theoretical properties and practical performance of these schemes are yet to be analyzed in greater depth.

The idea of correcting the mismatch of a knowledge-driven model with a data-driven model is akin to hybrid semi-parametric modeling \citep{Thompson1994,VonStosch2014}, specifically a parallel hybrid model structure. The consideration of non-parametric models, whereby the nature and number of parameters is not determined by a priori knowledge but tailored to the data at hand, makes perfect sense to capture the structural plant-model mismatch in RTO applications. In principle, this approach is even amenable to a completely model-free RTO scheme by simply discarding the mechanistic model component. But the effect of removing this mechanistic knowledge in a practical RTO setup has seldom been investigated to date.

\nomenclature[A]{GP}{Gaussian process}
\nomenclature[A]{MPC}{model predictive control}

Model predictive control (MPC) is closely related to RTO in that these two technologies entail the repeated solution of a model-based optimization problem at their core \citep{Rawlings2017}. Similar to RTO, a majority of successful MPC implementations have so far relied on mechanistic models. But there has been a renewal of interest in data-driven approaches, which use surrogate models trained on historical data or mechanistic model simulations to drive the optimization. The type of surrogate models used in MPC include  artificial neural networks \citep{Piche2000,Wu2019} and GPs \citep{Kocijan2004}. However, comparatively little work has been published on embedding hybrid models into MPC in order to reduce the dependency on data and infuse physical knowledge for better extrapolation capability \citep{Klimasauskas1998,Zhang2019}. 

A recent trend in MPC has been to include learning or self-reflective objectives alongside control performance objectives \citep{Hewing2020}. Self-reflective MPC seeks to minimize the controller's own performance loss in the presence of uncertainty \citep{Feng2018}. Instead, learning objectives aim to promote accurate future state and parameter estimates, inspired by optimal experiment design or persistent excitation ideas \citep{Larsson2013,Heirung2015,Marafioti2014}. In data-driven MPC for instance, recent research has investigated on-line learning of the surrogates to improve performance and reliability, with a particular interest in GPs \citep{Maiworm2018,Bradford2019,Bradford2020}. In essence, MPC with learning seeks to strike a balance between exploitation against exploration, which is akin to the dual control problem \citep{Wittenmark1995} and is also the central paradigm in the fast-developing field of reinforcement learning \citep{Speilberg2019,Kim2020,
petsagkourakis2020chance, Petsagkourakis2020}. Likewise, several modifier-adaptation schemes have incorporated excitation terms in the constraints of the RTO model in order to enable more accurate gradient estimates from noisy measurements \citep{Marchetti2010,Rodger2011}. But the vast potential of machine learning and reinforcement learning has remained largely untapped in the RTO context \citep{Powell2020}.

Other areas closely related to real-time optimization comprise black-box optimization and surrogate-based optimization, which find many applications in process flowsheeting, computational fluid dynamics, or molecular dynamics \citep{Biegler2014}. They can be broadly classified into local and global approaches. Global approaches proceed by constructing a surrogate model based on an ensemble of simulations before optimizing it, often within an iteration where the surrogate is progressively refined. A number of practical implementations rely on neural networks \citep{Henao2011}, GPs \citep{Caballero2008,Quirante2015,Kessler2019}, or a combination of various basis functions \citep{Wilson2017,Boukouvala2017} for the surrogate modeling. Bayesian optimization has gained significant popularity for tackling problems with expensive function evaluations, with prominent algorithms such as efficient global optimization \citep{Jones1998} and sequential kriging optimization \citep{Huang2006} that leverage GP surrogates and so-called acquisition functions to strike a balance between exploitation and exploration. Radial basis function (RBF) surrogates have also proven effective to optimize expensive black-box function \citep{Gutmann2001,Costa2018}. Handling constrained problems with this class of methods \added{still }constitutes an active field of research \deleted{nonetheless }\citep{Audet2018,Cartis2018}.

\nomenclature[A]{RBF}{radial basis function}

By contrast, local approaches \added{seek to} maintain an accurate approximation of the original optimization problem within a trust region, whose position and size are adapted iteratively. This procedure entails updating or reconstructing the surrogate model as the trust region moves around, but it benefits from a well-developed convergence theory providing sufficient conditions for local optimality in unconstrained and bound-constrained problems \citep{Conn2000,Conn2009,March2012,Cartis2019}. Extensions of these approaches to constrained flowsheet optimization include the work by \citet{Eason2016,Eason2018} and \citet{Bajaj2018}, while constrained multi-fidelity optimization was considered by \citet{March2012b}. In particular, the latter uses GP surrogates as low-fidelity models and their adaptation is akin to modifier adaptation with GP surrogates as developed by \citet{Ferreira2018} and \citet{delRio2019}. These connections between the modifier-adaptation and trust-region frameworks were also delineated in a short note by \citet{Bunin2014}. But while integrating local and global concepts from surrogate-based optimization methods within modifier adaptation is indeed appealing, this integration should account for the added complexity posed by noisy process data or changing optima over time in RTO. \added{Further developments in this area include probabilistic derivative-free trust-region methods (\citealp{Bandeira2014}; \citealp{Larson2016}; \citealp{Chen2018}), which rely on randomized surrogate models and can efficiently handle uncertainty. These ideas were recently connected to GP surrogates by \citet{augustin2017} and then modifier adaptation by \citet{Shukla2020}.}

Considering all this, the main focus of this paper is on improving modifier-adaptation schemes in terms of speed and reliability by integrating concepts and ideas from the areas of Bayesian optimization and derivative-free optimization. Specifically, the proposed modifier-adaptation schemes embed a physical model and trust-region concepts to minimize risk during the exploration, while relying on GPs to capture the plant-model mismatch in a non-parametric way and drive the exploration by means of acquisition functions. Key elements of novelty include the adaptation of the trust region based on the GPs' mean predictor ability to capture the plant-model mismatch in the cost and constraints and the exploitation of the GPs' variance estimators to maintain sufficient excitation during the search. \added{The focus is on algorithms that target good practical performance, rather than providing global convergence certificates at the cost of practicality. The performance of the proposed schemes is analyzed by means of numerical examples, including the benefits of using an acquisition function, knowing the process noise, or specifying a prior}\deleted{We furthermore investigate the effect of removing the prior,} knowledge-based model\deleted{ component and the effect of process noise by means of numerical examples}. 

The rest of the paper provides background on MA and GP in Section~\ref{sec:prelim}, then presents and analyses the new modifier-adaptation algorithm in Section~\ref{sec:meth}. This algorithm is illustrated with a simple quadratic optimization problem throughout Section~\ref{sec:meth} and with practical case studies in Section~\ref{sec:case}, before drawing final remarks in Section~\ref{sec:concl}.

\section{Preliminaries}
\label{sec:prelim}

\subsection{Modifier Adaptation}
\label{sec:MA}

The problem of optimizing the steady-state performance of a given plant subject to operational or safety constraints can be formulated as:
\begin{align}
\min_{{\bf u}\in\mathcal{U}}~ & G^{\rm p}_0\left({\bf u}\right) := g_0\left({\bf u},{\bf y}^{\rm p}({\bf u})\right)\label{eq:plant_problem}\\
\text{s.t.}~ & G_i^{\rm p}\left({\bf u}\right) := g_i\left({\bf u},{\bf y}^{\rm p}({\bf u})\right)\leq 0, \quad i=1\ldots n_g \nonumber
\end{align}
where ${\bf u}\in\mathbb{R}^{n_{u}}$ and ${\bf y}^{\rm p}\in\mathbb{R}^{n_{y}}$ are vectors of the plant input and output variables, respectively; $g_i:\mathbb{R}^{n_{u}}\times\mathbb{R}^{n_{y}}\rightarrow\mathbb{R}$, $i=0...,n_{g}$, denote the cost and inequality constraint functions; and $\mathcal{U}\subseteq\mathbb{R}^{n_u}$ is the control domain, e.g. lower and upper bounds on the input variables, ${\bf u}^{\rm L}\leq{\bf u}\leq{\bf u}^{\rm U}$. Notice the superscript $\left(\cdot\right)^{\rm p}$ used to indicate plant-related quantities.

The RTO challenge is of course that an exact mapping ${\bf y}^{\rm p}(\cdot)$ is unknown in practice, and the output ${\bf y}^{\rm p}({\bf u})$ can only be measured for a particular input value  ${\bf u}$, in the manner of a noisy oracle. However, provided that a\deleted{non-ideal (approximate)} model of the plant's input-output behavior is available, represented by the parametric function ${\bf y}({\bf u},\cdot)$, one may solve the following model-based optimization problem instead:
\begin{align}
\min_{{\bf u}\in\mathcal{U}}~ & G_0\left({\bf u}\right) := g_0\left({\bf u},{\bf y}({\bf u},\boldsymbol{\theta})\right)\label{eq:model_problem}\\
\text{s.t.}~ & G_i\left({\bf u}\right) := g_i\left({\bf u},{\bf y}({\bf u},\boldsymbol{\theta})\right)\leq 0, \quad i=1\ldots n_g \nonumber
\end{align}
where $\boldsymbol{\theta}\in\mathbb{R}^{n_\theta}$ is a vector of adjustable model parameters. 

In the presence of plant-model mismatch and process disturbances, the optimal solution value of Problem~\eqref{eq:model_problem} could be significantly different from that of Problem~\eqref{eq:plant_problem}. For this reason, a traditional two-step RTO scheme would try to reduce the plant-model mismatch by adjusting (a subset of) the model parameters with new plant measurements collected at each iteration. However, the convergence of such a scheme to a plant optimum is dependent upon a model adequacy condition \citep{Forbes1994,Chachuat2009}, whereby the model and plant optima match for at least one set of parameter values.

By contrast, the measurements in a modifier-adaptation scheme are used to correct the cost and constraint function values at a given iterate ${\bf u}^k$, in order to determine the next input or set-point values ${\bf u}^{k+1}$ \citep{Marchetti2009}:
\begin{align}
{\bf u}^{k+1} \in \argmin_{{\bf u}\in\mathcal{U}}~ & G_0({\bf u}) + (\boldsymbol{\lambda}_{\delta G_0}^k)^{\intercal} {\bf u} \label{eq:modified_problem}\\
\text{s.t.}~ & G_i\left({\bf u}\right) + \varepsilon_{\delta G_i}^k + (\boldsymbol{\lambda}_{\delta G_i}^k)^{\intercal} [{\bf u}-{\bf u}^k] \leq 0, \quad i=1\ldots n_g \nonumber
\end{align}
where $\varepsilon^{k}_{\delta G_i}\in\mathbb{R}$ are zeroth-order modifiers for the constraints, and $\boldsymbol{\lambda}_{\delta G_i}^k\in\mathbb{R}^{n_u}$ are first-order modifiers for the cost and constraints. The use of modifiers is appealing in that a KKT point ${\bf u}^\infty$ for the corrected model-based problem \eqref{eq:modified_problem} is also a KKT point for the original problem \eqref{eq:plant_problem}, provided that the modifiers satisfy \citep{Marchetti2009}:
\begin{align*}
  \varepsilon^{k}_{\delta G_i} =\ & G_i^{\rm p}({\bf u}^\infty) - G_i({\bf u}^\infty), \quad i=1\ldots n_g \\
  \boldsymbol{\lambda}_{\delta G_i}^k =\ & \boldsymbol{\nabla} G_i^{\rm p}({\bf u}^\infty) - \boldsymbol{\nabla} G_i({\bf u}^\infty), \quad i=0\ldots n_g
\end{align*}
A simple update rule for the modifiers that fulfills the foregoing conditions upon convergence is:
\begin{align*}
\varepsilon^{k+1}_{\delta G_i} =\ & (1-\eta) \varepsilon^{k+1}_{\delta G_i} + \eta\left[G_i^{\rm p}({\bf u}^k) - G_i({\bf u}^k)\right] \\
\boldsymbol{\lambda}_{\delta G_i}^{k+1} =\ & (1-\eta) \boldsymbol{\lambda}_{\delta G_i}^k + \eta\left[\boldsymbol{\nabla} G_i^{\rm p}({\bf u}^k) - \boldsymbol{\nabla} G_i({\bf u}^k)\right]
\end{align*}
where the tuning parameters $\eta\in(0,1]$ may be reduced to help stabilize the iterations. Apart from choosing a suitable $\eta$, the biggest burden with this approach is estimating the gradients $\boldsymbol{\nabla} G_i^{\rm p}({\bf u}^k)$ of the cost and constraint functions at each RTO iteration. A range of methods were reviewed in the paper's introduction to assist with this estimation. \added{Approaches to enforcing model adequacy in modifier-adaptation schemes are also available, for instance by means of a tailored parameter estimation procedure \citep{Ahmad2019}.}

\subsection{Gaussian Processes and Acquisition Functions}
\label{sec:GP}

GP regression is a method of interpolation developed by \citet{Krige1951} and popularized by the machine learning community \citep{Rasmussen2006}. It aims to describe an unknown function $f:\mathbb{R}^{n_u}\to\mathbb{R}$ using noisy observations, $y = f({\bf u}) + \nu$, where $\nu\sim\mathcal{N}(0,\sigma_{\nu}^{2})$ is Gaussian distributed measurement noise with zero mean and (possibly unknown) variance $\sigma_{\nu}^{2}$. GPs themselves consider a distribution over functions and may be regarded as a generalization of multivariate Gaussian distributions:
\begin{align*}
f(\cdot) \sim \mathcal{GP}(m(\cdot),k(\cdot,\cdot))
\end{align*}
where the mean function $m(\cdot)$ can be interpreted as the deterministic part of the function; and the covariance function $k(\cdot,\cdot)$ accounts for correlations between the function values at different points. 

One popular choice for the covariance function is the squared-exponential (SE) kernel \citep{Rasmussen2006}:
\begin{align*}
k({\bf u},{\bf u}') :=\ & \sigma_n^2\exp\left(-\frac{1}{2}({\bf u}-{\bf u}')^{\intercal}\boldsymbol{\Lambda}({\bf u}-{\bf u}')\right)
\end{align*}
where $\sigma_{n}^{2}$ is the covariance magnitude; and $\boldsymbol{\Lambda} := {\rm diag}(\lambda_1\cdots\lambda_{n_u})$ is a scaling matrix. Underlying this kernel choice is the assumption that the inferred function $f$ is both smooth and stationary. But other kernels could of course be selected, such as the Mat\'ern class of covariance functions \citep{Rasmussen2006}. We furthermore choose a constant mean function:
\begin{align*}
m({\bf u}) :=\ & c
\end{align*}
where $c$ is the scalar offset. This choice is motivated by the fact that since GPs are used to describe the plant-model mismatch in modifier adaptation, it is safe for their predictions to tend to a constant offset when extrapolating away from the measurement points\added{ \citep{Thompson1994}}.

Maximum likelihood estimation is commonly applied to infer a GP's hyperparameters $\boldsymbol{\Psi} := [c~\sigma_n~\sigma_{\nu}~\lambda_1\:\ldots\:\lambda_{n_u}]^{\intercal}$, where $\sigma_{\nu}$ may be excluded in case the measurement noise variance is known. Consider $N$ (noisy) function observations, denoted by ${\bf y} := [y_1\:\cdots\:y_{N})]^{\intercal} \in \mathbb{R}^{N}$, with corresponding inputs gathered in the matrix ${\bf U} := \left[ {\bf u}_1\:\cdots\:{\bf u}_{N} \right] \in \mathbb{R}^{n_u \times N}$. The log-likelihood of the observed data, ignoring constant terms, is given by:
\begin{align*}
\mathcal{L}(\boldsymbol{\Psi}) := -\frac{1}{2}\ln(|{\bf K({\bf U})}|) - \frac{1}{2}({\bf y}-{\bf 1}c)^{\intercal}\,{\bf K}({\bf U})^{-1}\,({\bf y}-{\bf 1}c)
\end{align*}
with $K_{ij}({\bf U}) := k({\bf u}_i,{\bf u}_j) + \sigma_{\nu}^{2}\delta_{ij}$ for all $(i,j) \in \{1\ldots N\}^2$; and Kronecker's delta function $\delta_{ij}$.

The predicted distribution of $f({\bf u})$ at an arbitrary input point ${\bf u}$, given the input-output data $({\bf U},{\bf y})$ and the maximum-likelihood estimates of $\boldsymbol{\Psi}$, follows a Gaussian distribution:
\begin{align}
f({\bf u}) \mid {\bf U},{\bf y} \sim \mathcal{N}(\mu_f({\bf u}),\sigma_f^2({\bf u}))
\label{eq:fdist}
\end{align}
where the posterior mean function $\mu_f$ and the posterior variance function $\sigma^2_f$ are computed as:
\begin{align*}
\mu_f({\bf u}) :=\ & {\bf r}({\bf u},{\bf U})\, {\bf K}({\bf U})^{-1}\, {\bf y} + c \\
\sigma_f^2({\bf u}) :=\ & \sigma_n^2 - {\bf r}({\bf u},{\bf U})\, {\bf K}({\bf U})^{-1}\, {\bf r}({\bf u},{\bf U})^{\intercal}
\end{align*}
with ${\bf r}({\bf u},{\bf U}) := [k({\bf u},{\bf u}_1)\:\cdots\:k({\bf u},{\bf u}_{N})]$.

\begin{figure}[tb]
\centering
\includegraphics[scale = 0.46]{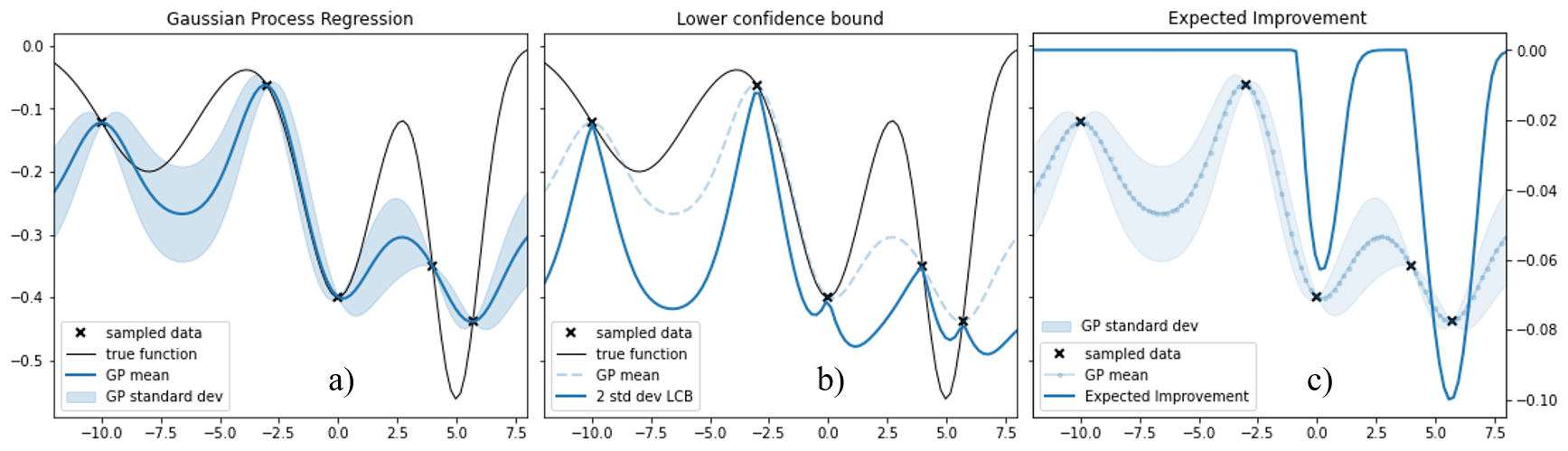}
\caption{({\bf a}) Illustration of the posterior mean and variance functions in the Gaussian process regression of a sampled function (Equation~\ref{eq:fdist}). ({\bf b}) Corresponding lower confidence bound (LCB) acquisition function (Equation~\ref{eq:LCB}\added{ with $\beta=3$}). ({\bf c}) Corresponding expected improvement (EI) acquisition function (Equation~\ref{eq:EI}).}
\label{fig:BayesOpt}
\end{figure}

In practice, the mean $\mu_f$ corresponds to the GP's prediction at ${\bf u}$, while the variance $\sigma_f^2$ provides a measure of the uncertainty associated to this prediction (Figure~\ref{fig:BayesOpt}a). Both functions are exploited in so-called acquisition functions, which constitute the workhorse of Bayesian estimation in balancing exploration versus exploitation \citep{Shahriari2016}. Two popular acquisition functions are reviewed next, namely {\em lower confidence bound} (LCB) and {\em expected improvement} (EI). Theses will be considered \deleted{later }as objective functions in the optimization subproblems of the modifier-adaptation algorithm\added{ (cf. Section~\ref{sec:meth})}. Alternative acquisition functions include {\em probability of improvement} \citep{Kushner1964}, {\em knowledge gradient} \citep{Frazier2009}, and {\em entropy search} \citep{Hennig2012}\deleted{ to name but a few}.\added{ The focus herein is on showing the benefits of using an acquisition function, but a more detailed comparison between various acquisition functions is left for future research.}

\nomenclature[A]{LCB}{lower confidence bound}
\nomenclature[A]{EI}{expected improvement}

\paragraph{Lower Confidence Bound}

With the notation introduced previously (Equation~\ref{eq:fdist}), this acquisition function is given by (Figure~\ref{fig:BayesOpt}b):
\begin{align}
\label{eq:LCB}
\mathcal{A}_{\rm LCB}[\mu_f,\sigma_f]({\bf u}) := \mu_f({\bf u}) - \beta \sigma_f({\bf u})
\end{align}
where $\beta$ may be interpreted as an exploration weight. Notice the negative sign of the exploration term in Equation~\eqref{eq:LCB}, which is consistent with the formulation of the RTO Problem~\eqref{eq:plant_problem} as a minimization. The LCB function is based on the principle of optimism in the face of uncertainty, with a view to minimizing regret. Its early use can be traced back to the work by \citet{Lai1985} on rule allocations, and later by \citet{Agrawal1995} in the context of reinforcement learning. \added{An important appeal of this acquisition function lies in its simplicity.}

\paragraph{Expected Improvement}

This acquisition function is expressed as (Figure~\ref{fig:BayesOpt}c):
\begin{align}
\label{eq:EI}
\mathcal{A}_{\rm EI}[\mu_f,\sigma_f,f_{\rm L}]({\bf u}) := - \left[f_{\rm L}-\mu_f({\bf u})\right] \Phi\left( \frac{f_{\rm L}-\mu_f({\bf u})}{\sigma_f({\bf u})} \right) - \sigma_f({\bf u})\ \phi\left( \frac{f_{\rm L}-\mu_f({\bf u})}{\sigma_f({\bf u})} \right)
\end{align}
where $\phi(\cdot)$ and $\Phi(\cdot)$ are the standard normal probability density and cumulative distribution functions, respectively; and $f_{\rm L} := \min (y_1,\ldots,y_N)$ is the best observed value, possibly replaced with the lowest mean value, $\min (\mu_f({\bf u}_1),\ldots,\mu_f({\bf u}_N))$, in case the observations carry significant noise. This expression corresponds to $\mathbb{E}[\max(f_{\rm L}-\mu_f(\cdot),0)]$, where the improvement function $\max(f_{\rm L}-\mu_f(\cdot),0)$ is only positive at points where the predicted mean value is lower than $f_{\rm L}$. The negative signs are introduced so that $\mathcal{A}_{\rm EI}$ can be used as objective function in a minimization problem\deleted{, rather than maximized}. Its introduction is credited to \citet{Mockus1975} and it was later popularized via the efficient global optimization (EGO) algorithm by \citet{Jones1998}. 
\bigskip

Both the LCB and EI acquisition functions seek to balance exploration and exploitation in order to reduce the overall number of observations. Computational benchmarks tend to favor EI over LCB though, since the latter may lead to excessive exploration \citep{Snoek2012,Shahriari2016}. Nevertheless, there are theoretically motivated guidelines for tuning the weight $\beta$ to achieve optimal regret \citep{Srinivas2010}, and thereby boost the performance of \replaced{LCB}{USB}. \replaced{Further practical caveats are that}{Both functions furthermore come with caveats in practice, as} LCB typically comprises a larger number of local optima, whereas EI can present large flat areas. These characteristics call for randomized search or complete search approaches in applications \citep[see, e.g.][]{Torn1989,Schweidtmann2020}.

\section{Methodology}
\label{sec:meth}

\subsection{Modifier-Adaptation Algorithm Statement}
\label{sec:MA-GP-ITR}

The use of GPs to describe the plant-model mismatch in an RTO problem was first proposed by \citet{Ferreira2018}. The main idea \replaced{is}{was} for these GP modifiers to correct the cost and each constraint separately:
\begin{align*}
G_i^{\rm p}-G_i \sim \mathcal{GP}\left(\mu_{\delta G_i},\sigma_{\delta G_i}^{2}\right), \quad i=0\ldots n_g
\end{align*}
The\added{n, the} following modified optimization problem \replaced{is}{was then} solved in an RTO iteration:
\begin{align}
{\bf u}^{k+1} \in \argmin_{{\bf u}\in\mathcal{U}}~ & [G_0+\mu^k_{\delta G_0}]({\bf u}) \label{eq:modified_problem_GP}\\
\text{s.t.}~ & [G_i+\mu^k_{\delta G_i}]({\bf u}) \leq 0,\quad i=1\ldots n_g \nonumber
\end{align}
where $\mu^k_{\delta G_i}$ denotes the mean of the GP trained with the input-output data set $({\bf U}^k,\boldsymbol{\delta}{\bf G}_i^k)$; and $\boldsymbol{\delta}{\bf G}_i^k$ comprises measurements of the mismatch $\delta G_i(\cdot):= G^{\rm p}_i(\cdot)-G_i(\cdot)$ for inputs in the matrix ${\bf U}^k$. \added{In addition to providing zeroth- and first-order correction terms for the $G_i$'s, observe that GP modifiers are also capable of second-order corrections. Such corrections of the curvature of a prior model are appealing insofar as they can help enforce the model adequacy conditions of modifier adaptation (\citealp{Marchetti2009}; \citealp{Gottumukkula2020}), that is, positive semi-definiteness of the reduced Hessian of model-based optimization problem at a plant optimum.} This idea of \replaced{using GPs to construct global, nonlinear}{constructing GP} surrogates for \deleted{the cost and constraint} black-box functions is also shared by various derivative-free algorithms (\citealp{March2012b}; \citealp{Picheny2016}\added{; \citealp{augustin2017}}).

Herein, we revisit this idea by introducing trust-region concepts from the fields of derivative-free and surrogate-based optimization together with acquisition functions from Bayesian optimization. The modified optimization problem that is solved in each RTO iteration becomes:
\begin{align}
{\bf d}^{k+1} \in \argmin_{\bf d}~ & \mathcal{A}[G_0+\mu^k_{\delta G_0},\sigma^k_{\delta G_0},\cdot]({\bf u}^k+{\bf d}) \label{eq:modified_problem_GP+TR}\\
\text{s.t.}~ & [G_i+\mu^k_{\delta G_i}]({\bf u}^k+{\bf d}) \leq 0,\quad i=1\ldots n_g \nonumber\\
& \|{\bf d}\|\leq\Delta^{k}, \quad {\bf u}^k+{\bf d}\in\mathcal{U} \nonumber
\end{align}
where $\Delta^k\added{>} 0$ is the trust-region radius for the predicted step ${\bf d}^{k+1}\in\mathbb{R}^{n_u}$; and $\mathcal{A}$ is an acquisition function for the cost predictor $G_0+\mu^k_{\delta G_0}$ and the associated error estimate $\sigma^k_{\delta G_0}$\replaced{. Subsequently, this acquisition function will be chosen as}{, which may be} either the LCB or EI function (cf. Section~\ref{sec:GP}) or the cost predictor itself if exploration is not considered.

Solving Problem~\eqref{eq:modified_problem_GP+TR} is akin to conducting a constrained Bayesian optimization within a trust-region. The various steps used to adapt this trust region and handle the constraints are summarized in Algorithm~\ref{alg:MA-GP-TR} and commented below.

\begin{algorithm}[tbp]
\caption{Modifier adaptation with Gaussian process, trust region and acquisition function}\label{alg:MA-GP-TR}
\TabPositions{8cm}
{\bf Input:} initial data sets $({\bf U}^0,\boldsymbol{\delta}{\bf G}_i^0)$, $i=0\ldots n_g$; trained GP modifiers $\mu^0_{\delta G_i}$, $i=0\ldots n_g$ and $\sigma^0_{\delta G_0}$; initial operating point $\mathbf{u}^0\in\mathcal{U}$; initial and maximal trust-region radii $0<\Delta^0<\Delta_{\rm max}$; trust-region parameters $0<\eta_1<\eta_2<1$, $0 < \gamma_{\rm red} < 1 < \gamma_{\rm inc}$, and $\added{\alpha}>0$; subset of unrelaxable constraints $\mathcal{UC}\subseteq \{1\ldots n_g\}$\\

{\bf Repeat: for $k=0,1,\ldots$} 
\begin{enumerate}
\item Check criticality\\
If $\Delta^{k} > \added{\alpha} \left\| \boldsymbol{\nabla}_{\rm red} [G_0+\mu^k_{\delta G_0}]({\bf u}^k) \right\|$: \quad $\Delta^{k}\ \leftarrow\  \gamma_{\rm red} \Delta^{k}$

\item Solve modified optimization problem (Equation~\ref{eq:modified_problem_GP+TR})
\tab $\rhd$\ $\mathbf{d}^{k+1}$

\item Get process cost and constraint measurements
\tab $\rhd$\ $G^{\rm p}_i({\bf u}^{k}+{\bf d}^{k+1}), i=0\ldots n_g$

\item Check infeasibility\\
If \added{either }Problem~\eqref{eq:modified_problem_GP+TR} is infeasible, or \added{if }$G_{i}^{\rm p}({\bf u}^{k}+{\bf d}^{k+1})>0$ for any $i\in\mathcal{UC}$:\\
$\hphantom{\quad} \Delta^{k+1}\ \leftarrow\ [\gamma_{\rm red},1] \Delta^{k}$, \quad ${\bf u}^{k+1}\ \leftarrow\ {\bf u}^k\ \text{(reject)}$, \quad and goto Step 7

\item Compute merit function  (Equation~\ref{eq:merit_function})
\tab $\rhd$\ $\rho^{k+1}$

\item Update trust region\vspace{-.5em}
\begin{align*}
& \text{If $\rho^{k+1} > \eta_2\ \wedge\ \|{\bf d}^{k+1}\|=\Delta^{k}$:} && \Delta^{k+1}\ \leftarrow\ \gamma_{\rm inc} \Delta^{k}, && {\bf u}^{k+1}\ \leftarrow\ {\bf u}^k + {\bf d}^{k+1}\ \text{(accept)} \qquad\qquad\quad\\
& \text{Else If $\rho^{k+1} < \eta_1$:} && \Delta^{k+1}\ \leftarrow\ \gamma_{\rm red} \Delta^{k}, && {\bf u}^{k+1}\ \leftarrow\ {\bf u}^k\ \text{(reject)}\\
& \text{Else:} && \Delta^{k+1}\ \leftarrow\  \Delta^{k}, && {\bf u}^{k+1}\ \leftarrow\ {\bf u}^k + {\bf d}^{k+1}\ \text{(accept)}
\end{align*}

\item Update data sets
\tab $\rhd$\ $({\bf U}^{k+1},\boldsymbol{\delta}{\bf G}_i^{k+1}), i=0\ldots n_g$
\item Update GP modifiers
\tab $\rhd$\ $\mu^{k+1}_{\delta G_i}$, $i=0\ldots n_g; \sigma^{k+1}_{\delta G_0}$
\end{enumerate}
\end{algorithm}

\paragraph{Initialization}

A set of GPs are trained on cost and constraint mismatch data in the initial step. There is considerable freedom regarding the choice of this initial training set $({\bf U}^0,\boldsymbol{\delta}{\bf G}_i^0)$, $i=0\ldots n_g$ as well as the initial trust-region center ${\bf u}^0$ and radius $\Delta^0$. One approach entails defining the initial trust region first, then selecting an initial sample set within this trust region in a second step. Such an initial trust region may leverage process knowledge and physical insight in practice. Identifying a feasible starting point for (a subset of) the process constraints could also be via the solution of an auxiliary feasibility problem prior to running Algorithm~\ref{alg:MA-GP-TR} \citep{Bajaj2018}. Sample points may then be generated within this trust region by imposing finite perturbations along each input direction or using quasi-random sampling, ideally so that the GP surrogates can be certified to be \added{(probabilistically)} fully linear---further discussions of the full linearity property are deferred to the convergence subsection below\added{ as well as the appendix}. When a process data set $({\bf U}^0,\boldsymbol{\delta}{\bf G}_i^0)$, $i=0\ldots n_g$ is preexisting, such as historical data, another approach involves constructing a maximal trust region that lies within a confidence percentile of the cost and constraint GP predictors $(\mu^0_{\delta G_i},\sigma^0_{\delta G_i})$ trained on this data set. Although such maximization problems are generally hard to solve because of their nonconvexity, good feasible solutions may nevertheless be obtained for practical purposes with any local solver and a multistart heuristic or using any global solver as feasibility pump \citep{Schweidtmann2020}.

The issue of scaling is closely related to that of trust-region and GP initialization. In practice, one can exploit the input domain $\mathcal{U}$ to scale the input variable to within $[0,1]$. The benefits of operating within a scaled input domain, both in terms of trust-region adaptation and GP training, are clear. A maximal trust-region radius may also be defined more conveniently in a scaled input domain, e.g. $\Delta_{\rm max}=0.7$. Note that there is furthermore considerable flexibility in the choice of the trust region parameters $\eta_1$, $\eta_2$, $\gamma_{\rm red}$ and $\gamma_{\rm inc}$. A common setting in trust-region methods, which is also the setting used for the numerical case studies below, is $\eta_1=0.2$, $\eta_2=0.8$, $\gamma_{\rm red}=0.8$ and $\gamma_{\rm inc}=1.2$. By contrast, the criticality parameter $\mu$ is problem dependent and may be set to an arbitrary large value if shrinking of the trust region upon convergence to a stationary point is not desirable.

\paragraph{Adaptation Mechanisms}

The trust region serves the dual purpose of restricting the step size to the neighborhood where the cost and constraint surrogates are deemed to be predictive, while also defining the neighborhood in which \replaced{additional}{the} points are sampled for \replaced{updating}{the construction of} these surrogates. The trust region update corresponds to Steps 1, 4 and 6 of Algorithm~\ref{alg:MA-GP-TR}. The latter comprises the classical update rules in trust-region algorithms \citep{Conn2009}, which is based on the ratio of actual cost reduction to predicted cost reduction:
\begin{align}
\rho^{k+1} := \frac{G_0^{\rm p}\left({\bf u}^{k}\right)-G_0^{\rm p}\left({\bf u}^{k}+{\bf d}^{k+1}\right)}{[G_0+\mu^k_{\delta G_0}]({\bf u}^{k})-[G_0+\mu^k_{\delta G_0}]({\bf u}^{k}+{\bf d}^{k+1})} \label{eq:merit_function}
\end{align}
The trust-region radius $\Delta^{k+1}$ is reduced whenever the accuracy ratio $\rho^{k+1}$ is too low. Conversely, $\Delta^{k+1}$ is increased if the optimization model \eqref{eq:modified_problem_GP+TR} takes a full step and the modified cost is deemed a good enough prediction of the plant cost variation around this point. Otherwise, the trust-region radius stays unchanged. As for the operating point update, the full step ${\bf d}^{k+1}$ is accepted when the accuracy ratio $\rho^{k+1}$ is large enough. Otherwise, the operating point remains unchanged, which would entail a back-tracking \added{from ${\bf u}^{k}+{\bf d}^{k+1}$ to ${\bf u}^{k}$} in a practical RTO setup. 

Before applying these updates, Step 4 asserts the feasibility of the modified optimization model \eqref{eq:modified_problem_GP+TR}\replaced{, and it stays at the same point and possibly reduces the trust region if this model is infeasible. Otherwise, it test the feasibility}{ and} of the plant constraints\replaced{ at the next point ${\bf u}^{k}+{\bf d}^{k+1}$, and rejects}{. Any infeasibility triggers a rejection of} the step ${\bf d}^{k+1}$ \replaced{with}{and} a possible reduction of the trust-region radius\replaced{ in case of infeasibility}{, again resulting in backtracking to point ${\bf u}^{k}$ in a practical RTO setup}. \deleted{Note that such a strategy also requires that the initial point ${\bf u}^0$ should satisfy all the plant constraints.} Backtracking is equivalent to the extreme-barrier approach in the trust-region literature \citep{Audet2006,Conn2009,Larson2019}, which assigns an infinite cost to points that violate any constraint. \added{Note that such a strategy also requires that the initial point ${\bf u}^0$ should satisfy all the plant constraints.} It is customary in this literature to distinguish between relaxable and unrelaxable constraints, where only the former may be violated along the search path. Various approaches to handling relaxable constraints within trust-region algorithms have been developed in recent years, including progressive-barrier, augmented-Lagrangian and filter methods \citep[see, e.g.][]{Picheny2016,Larson2019}. Integrating these techniques within a modifier-adaptation scheme is promising, but falls beyond the scope of the present paper. Instead, \replaced{Algorithm~\ref{alg:MA-GP-TR} applies}{we shall only apply} backtracking to the unrelaxable constraints\deleted{ subsequently} (subset $\mathcal{UC}$), while bypassing this check for the relaxable \deleted{inequality or equality }constraints.

The criticality test in Step 1 is inspired by state-of-the-art trust-region algorithms in derivative-free optimization. The aim is to keep the radius of the trust region comparable to some measure of stationarity in order for the surrogate model to become more accurate as the iterates get closer to a stationary point. The update of the trust-region radius in Step 1 forces it to converge to zero, hence defining a natural stopping criterion for this class of methods \citep{Conn2009}. In the presence of constraints, stationarity of the cost function may be substituted by Lagrangian stationarity or, alternatively, a reduced-gradient condition with:
\begin{align}
\label{eq:reduced_gradient}
\boldsymbol{\nabla}_{\rm red} [G_0+\mu^k_{\delta G_0}]({\bf u}^k) :=&\ \boldsymbol{\nabla} [G_0+\mu^k_{\delta G_0}]({\bf u}^k)\ {\bf N}^k
\end{align}
where the columns of ${\bf N}^k\in\mathbb{R}^{n_u\times(n_u-n_{g,{\rm a}})}$ form an orthogonal basis of the nullspace of the active constraint gradients at ${\bf u}^k$. However, it is better to treat this step as optional in a practical RTO setup, e.g. by allowing $\mu\to\infty$. This is because convergence of the trust-region radius to zero might hinder an RTO system's capability to react to process disturbances in order to track a time-varying optimum. \added{This is also the reason why the main iteration loop in Algorithm~\ref{alg:MA-GP-TR} does not specify a termination criterion.}

Apart from updating the trust region, both the data sets and the GP modifiers are updated at Steps 7 and 8, irrespective of whether the step ${\bf d}^{k+1}$ is accepted or not. The default strategy herein is to keep all of the past iterates and reconstruct the GPs by fitting all of their respective hyperparameters\added{, in the manner of a global surrogate model}. In order to prevent overfitting and numerical difficulties in constructing the GPs, \citet{Ferreira2018} proposed to keep a limited number of historical records in the input-output data set. This subset could comprise the $N$ most recent iterates or the $N$ nearest-neighbors to the next operating point ${\bf u}^{k+1}$. The former is akin to a forgetting strategy that is suitable for the tracking of a changing optimum, while the latter might be more appropriate to precisely locate a steady optimum. Moreover, the current iterate ${\bf u}^{k+1}$ need not be included in ${\bf U}^{k+1}$ should it be within a given radius of an existing point in ${\bf U}^{k}$, or ${\bf u}^{k+1}$ could be substituted for an existing nearby point in ${\bf U}^{k+1}$ instead. The computational burden of reconstructing the GPs at each iteration could furthermore be eased upon updating the covariance matrix at certain iterations only \citep{Rasmussen2006}. Another key RTO design decisions is whether to identify the measurement noise variance $\sigma_\nu^2$ alongside the other GP hyperparameters (cf. Section~\ref{sec:GP}), or to use an a priori noise variance provided by the sensor manufacturer or estimated from historical data. This discussion is deferred until the numerical analysis in Section~\ref{sec:analysis}.

\paragraph{Convergence and Performance Aspects}

\added{Derivative-free trust-region algorithms can be broadly classified into two categories, those with certified convergence and those which target good practical performance \citep{Conn2009}. The former are well established for unconstrained optimization problems and rely on the key property of (probabilistically) fully linear surrogate models---a summary of these results is reported in Appendix~\ref{app:convergence} for completeness. Global convergence is also certifiable for constrained optimization problems with such algorithms, e.g. by using penalty functions \citep{Larson2019}. However, the penalty approach cannot guarantee feasible iterates along the path to a critical point and it is therefore unsuitable for unrelaxable constraints in the RTO context. Another drawback of penalty methods is the need to update the penalty parameters, often by means of an outer loop, which can increase the number of function evaluations. Such trade-offs between convergence and performance are particularly relevant in the RTO context, where the optimum may change due to process disturbances or other external factors and progressing towards a process optimum sufficiently fast may be critical.}

\deleted{The modifier-adaptation scheme in Algorithm~\ref{alg:MA-GP-TR} is inspired by the derivative-free trust-region method in \citet[][Algorithm~10.1]{Conn2009}. The benefit of this design is that Algorithm~\ref{alg:MA-GP-TR} is globally convergent to a first-order critical point for unconstrained problems, upon imposing additional conditions such as full linearity of the surrogates in each trust-region subproblem and in the absence of noise; this convergence analysis is reported in Appendix~\ref{app:convergence} for completeness. Derivative-free trust-region methods that are provably convergent for problems with black-box constraints have also been developed in recent years \citep{Augustin2014,Echebest2017,Audet2018}. These methods are based on full-linearity assumptions as well.}

\deleted{Derivative-free trust-region methods can be broadly classified into two categories, those which target good practical performance and those for which convergence is established \citep{Conn2009}. The latter typically pay the price of practicality since ensuring full-linearity often requires taking extra sample points within the active trust region. This trade-off between convergence and performance is particularly relevant in the RTO context, whereby moving towards a process optimum sufficiently fast can be critical. With this in mind and the fact that the optimum can change due to process disturbances we do not enforce full linearity of the GP surrogates in Algorithm~\ref{alg:MA-GP-TR} and implement a simple backtracking strategy to handle constraint violation instead. Both the illustrative example in Section~\ref{sec:analysis} and the numerical case study in Section~\ref{sec:case} below confirm that Algorithm~\ref{alg:MA-GP-TR} can locate constrained process optima both efficiently and reliably, even in the presence of noise.}

\added{With this in mind, Algorithm~\ref{alg:MA-GP-TR} does not enforce (probabilistic) full linearity of the GP surrogates and it handles constraints through a simple backtracking (extreme barrier) approach rather than a penalty function. Consequently, the convergence of Algorithm~\ref{alg:MA-GP-TR} cannot be certified in general. But the algorithm steps are inspired by convergent derivative-free trust-region schemes (e.g., Algorithm~10.1 in \citealp{Conn2009}), so extra conditions such as probabilistic full linearity could be added in order to certify convergence to first-order critical points, at least in unconstrained RTO problems (cf. Appendix~\ref{app:convergence}).}

\deleted{Without enforcing full linearity of the surrogates, however, a derivative-free trust-region algorithm may get trapped around a suboptimal point.} We show with an illustrative example below (cf. Figures~\ref{fig:ex1_expl_no} \& \ref{fig:ex1_expl_obj}) that Algorithm~\ref{alg:MA-GP-TR} may indeed fail to steer the iterates to a process optimum when the cost and constraints are simply corrected in the manner of Problem~\eqref{eq:modified_problem_GP} in Step 2. A similar situation is known to occur in modifier-adaptation schemes that exploit past operating points in recursive gradient updates~\citep{Marchetti2010,Rodger2011}, where the addition of extra constraints in the RTO model to generate excitation can help mitigate the problem. Herein, we address this problem by leveraging ideas from Bayesian optimization for the first time. In particular, we use an acquisition function in order to promote exploration within the trust region (Problem~\ref{eq:modified_problem_GP+TR}). \deleted{The acquisition functions of interest, lower confidence bound (LCB) and expected improvement (EI), were described in Section \ref{sec:GP} and are assessed on an illustrative example in Section~\ref{sec:analysis}}. \added{Both the illustrative example in Section~\ref{sec:analysis} and the numerical case study in Section~\ref{sec:case} below confirm that Algorithm~\ref{alg:MA-GP-TR} can locate constrained process optima both efficiently and reliably with either the LCB or EI acquisition function (cf. Section \ref{sec:GP}).} \deleted{Acquisition functions could also be considered for the constraints of the RTO model themselves, either to prevent constraint violation or to improve the accuracy of the constraint GP surrogates \citep{Picheny2016,delRio2019}, although this falls beyond the scope of the present paper.} 

\added{Notice that acquisition functions could also be considered for the constraints of the modified optimization model \eqref{eq:modified_problem_GP+TR}. For instance, an LCB function could be used to relax the modified constraints:}
\begin{align*}
\added{[G_i+\mu^k_{\delta G_i}-\beta \sigma^k_{\delta G_i}]({\bf u}^k+{\bf d}) \leq 0,\quad i=1\ldots n_g}
\end{align*}
\added{thereby promoting further exploration within the trust region \citep{Picheny2016}. Conversely, tightening the modified constraints in the manner of an upper confidence bound (UCB):}
\begin{align*}
\added{[G_i+\mu^k_{\delta G_i}+\beta \sigma^k_{\delta G_i}]({\bf u}^k+{\bf d}) \leq 0,\quad i=1\ldots n_g}
\end{align*}
\added{could prevent constraint violations and thus reduce the need for backtracking during the RTO iterations \citep{delRio2019}. These variants are beyond the scope of the present paper and will be investigated as part of future work.}

\paragraph{Computational Aspects}

Traditional RTO systems often comprise complex numerical optimization subproblems as they rely on mechanistic models to drive the optimization. Correcting the cost and constraint functions with GP modifiers as in Problem~\eqref{eq:modified_problem_GP+TR} can introduce further nonlinearity and nonconvexity, thereby adding even more to this complexity. It is well known in particular that both the LCB and EI acquisition functions can exhibit a large number of local optima (cf.~Figure~\ref{fig:BayesOpt}). \replaced{Complete search methods in}{The application of} global optimization \replaced{are}{methods is} only computationally tractable for small-scale problems in practice\deleted{ \citep{Schweidtmann2020}}. Instead, the numerical case studies throughout this paper are solved using a local solver in combination with a multistart heuristic. The corresponding python codes are made available in the Supporting Information for the sake of reproducibility.

In principle, one could also decide to construct the GP modifiers from scratch, that is, without correcting an a priori mechanistic model. The optimization subproblems in such a model-free RTO system could be solved to guaranteed global optimality more efficiently using state-of-the-art complete-search algorithms \citep{Schweidtmann2020}. But the lack of a mechanistic model embedded into the optimization problem might significantly slow down the progress of the iterates to a plant optimum or be detrimental to the reliability of the RTO system. This trade-off is analyzed in greater details in the following section and later illustrated on the case studies\deleted{ too}.

\subsection{Algorithm Performance and Analysis: Illustrative Example}
\label{sec:analysis}

We \replaced{consider}{now analyze the performance of Algorithm~\ref{alg:MA-GP-TR} for various design choices by considering} the following simple optimization problem:
\begin{align}\label{eq:Illu1}
\min_{{\bf u}\in[-2,2]^2}\ \ & y_1({\bf u})\\
\text{s.t.}\ \ & y_2({\bf u}) \leq 0 \nonumber\\
& y_1({\bf u}) := u_1^2 + u_2^2 + \theta_1 u_1 u_2\nonumber\\
& y_2({\bf u}) := 1 - u_1 + u_2^2 + \theta_2 u_2 \nonumber
\end{align}
The (unknown) plant parameter values are taken as $\boldsymbol{\theta}^{\rm p} = [1\ 2]^{\intercal}$. The corresponding plant optimum (and the only KKT point here) is ${\bf u}^\ast \approx [0.368\ -0.393]^{\intercal}$, where the inequality constraint is active and the optimal cost is $y_1^\ast \approx 0.145$. In order to conduct the RTO, we assume that both outputs $y_1^{\rm p}$ and $y_2^{\rm p}$ are measured \replaced{but corrupted with}{and we add} a Gaussian white noise of variance $\sigma^2_{y_1}=\sigma^2_{y_2}=10^{-3}$\deleted{ to the simulated values}. Unless otherwise noted, we assume that the level of noise is not known a priori and therefore the variances $\sigma^2_{y_1}$ and $\sigma^2_{y_2} $ need to be estimated alongside the other GP hyperparameters (cf. Section~\ref{sec:GP}). We furthermore consider a nominal model with parameter values $\boldsymbol{\theta} = [0\ 0]^{\intercal}$, so that the problem presents a structural mismatch.

\begin{figure}[tb]
\begin{subfigure}{0.5\textwidth}
  \centering
  \includegraphics[width=.98\linewidth]{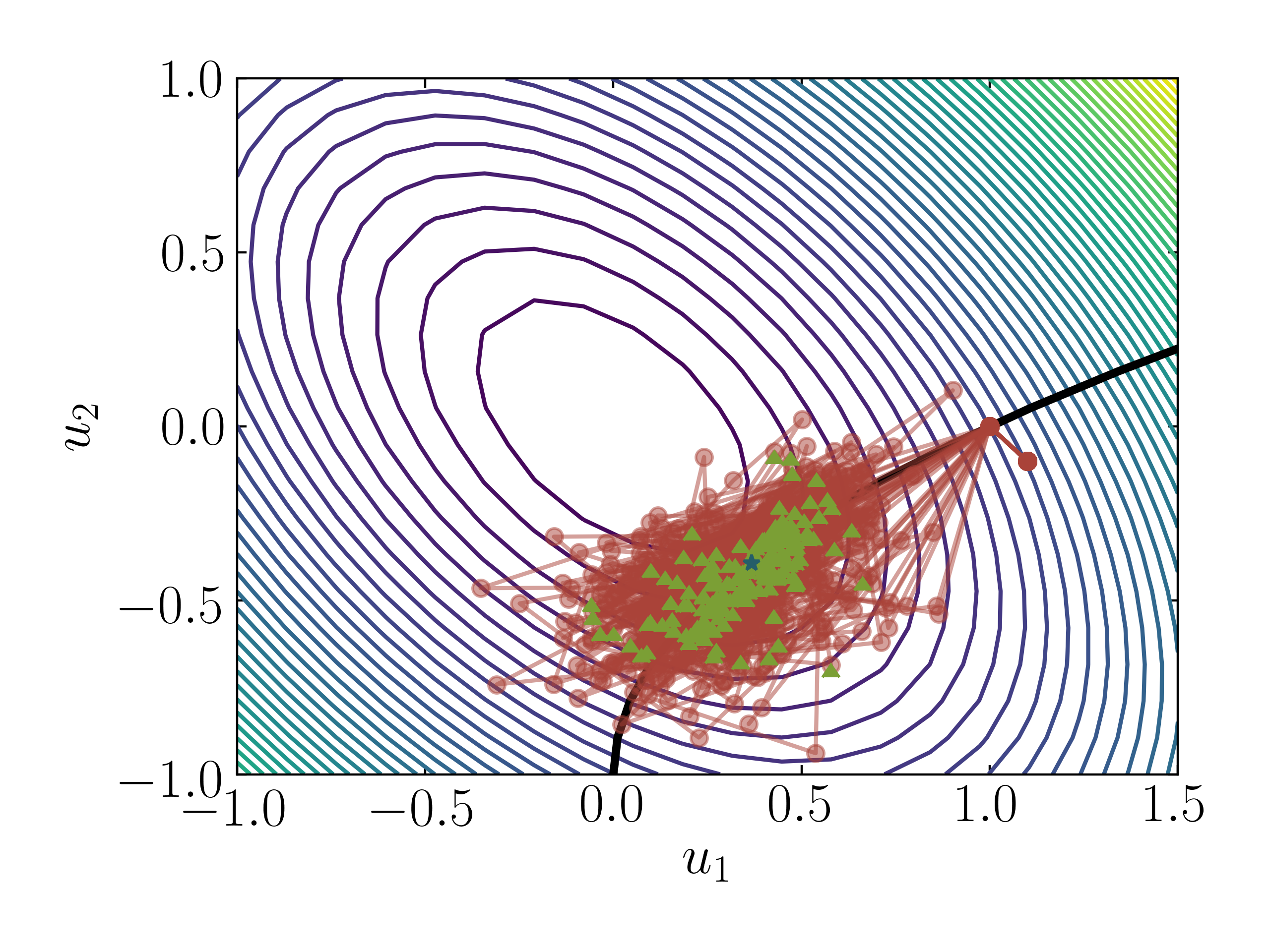}\vspace{-.85em}
  \caption{RTO iterations with $\eta=0.8$}
  \label{fig:ex1_std_08}
\end{subfigure}
\begin{subfigure}{0.5\textwidth}
  \centering
  \includegraphics[width=.98\linewidth]{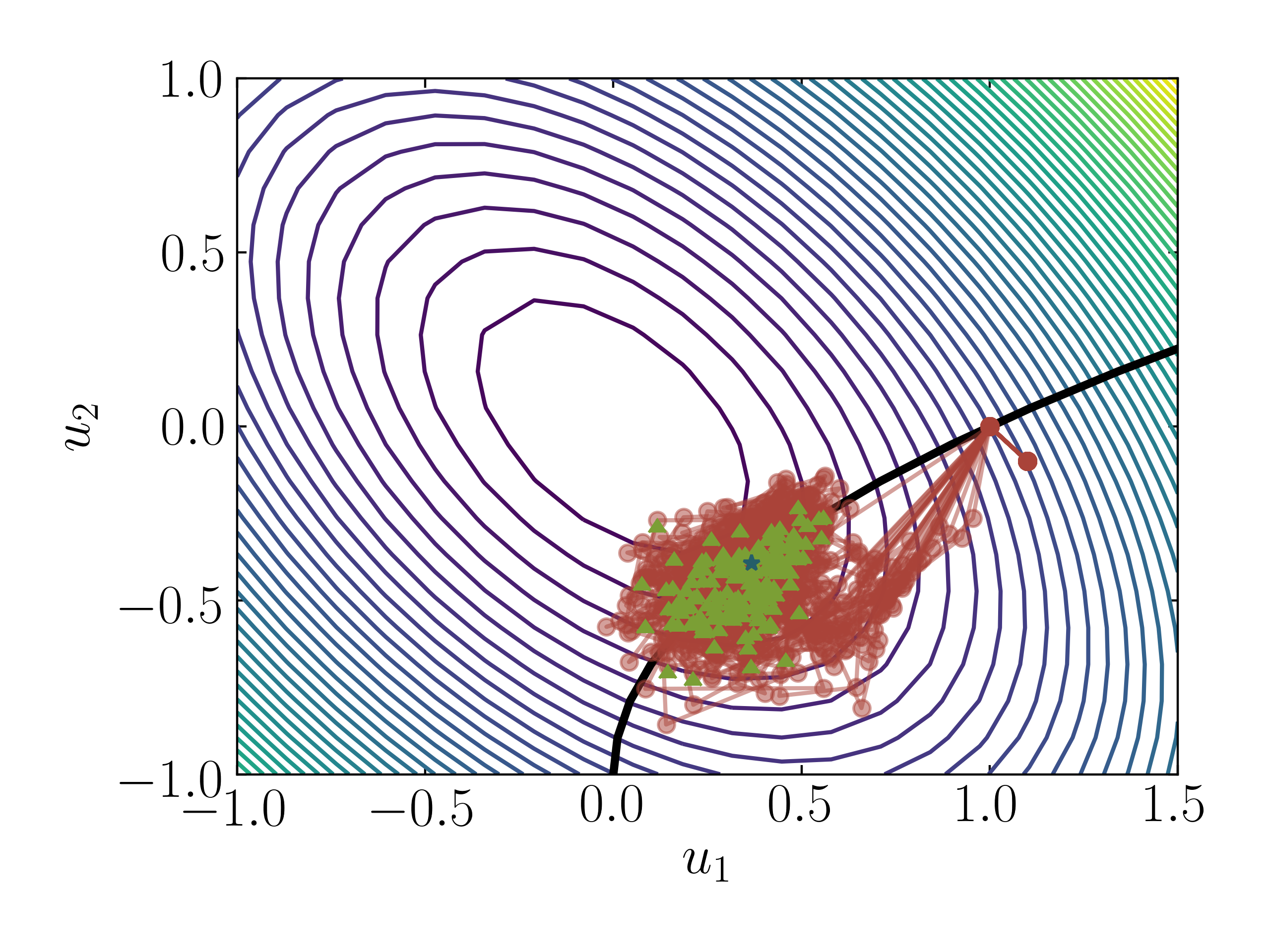}\vspace{-.85em}
  \caption{RTO iterations with $\eta=0.5$}
  \label{fig:ex1_std_05}
\end{subfigure}
\begin{subfigure}{0.5\textwidth}
  \centering
  \includegraphics[width=.98\linewidth]{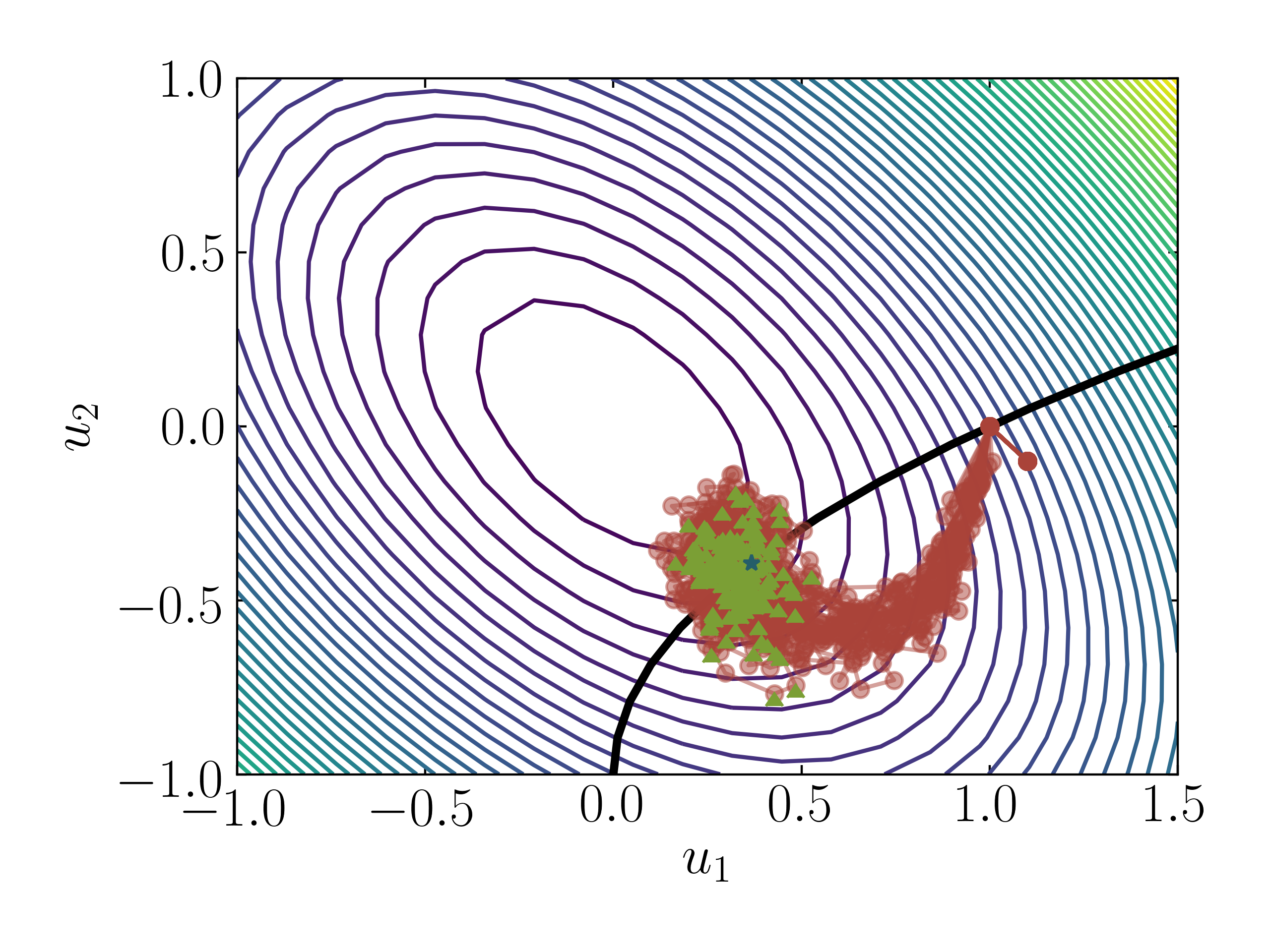}\vspace{-.85em}
  \caption{RTO iterations with $\eta=0.2$}
  \label{fig:ex1_std_02}
\end{subfigure}
\begin{subfigure}{0.5\textwidth}
  \centering
  \includegraphics[width=.98\linewidth]{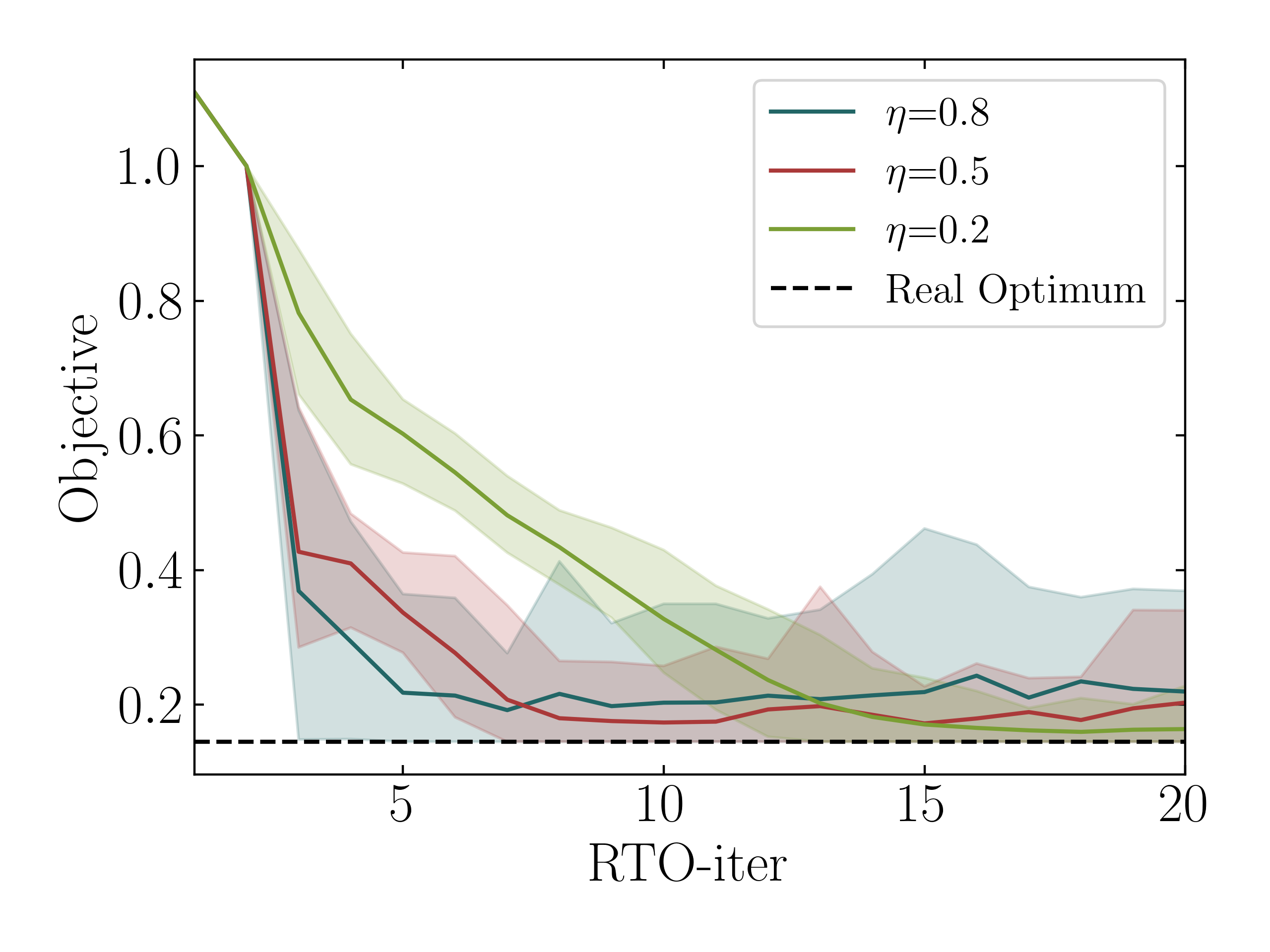}\vspace{-.85em}
  \caption{Comparative evolution of process cost with different gains}
  \label{fig:ex1_std_obj}
\end{subfigure}
\caption{\added{RTO iterations for Problem~\eqref{eq:Illu1} with a standard modifier-adaptation scheme based on the optimization model \eqref{eq:Illu1_ma0} with the update rules (\ref{eq:MA_filter_bias},\ref{eq:MA_filter_grad}) and finite-difference gradients. {\bf (a), (b), (c)} Clouds of iterates (red connected circles) for 30 process noise realizations, initialized from ${\bf u}^0 = [1.1\ -0.1]^{\intercal}$ and interrupted after 20 RTO iterations (green triangles); the process optimum is depicted with a blue star. {\bf (d)}~Evolution of the 95th percentile of process cost values over all the noise realizations with the RTO iterations (showing only the feasible iterates).} 
\label{fig:ex1_standard}}
\end{figure}

\added{In a traditional modifier-adaption scheme (cf. Section~\ref{sec:MA}), the following optimization model is solved in each RTO iteration:}
\begin{align}
\added{{\bf u}^{k+1} \in \argmin_{\bf u}~} & \added{y_1({\bf u})+(\boldsymbol{\lambda}^k_{\delta y_1})^\intercal {\bf u}} \label{eq:Illu1_ma0}\\
\added{\text{s.t.}\ \ } & \added{y_2({\bf u})+\varepsilon^k_{\delta y_2} +(\boldsymbol{\lambda}^k_{\delta y_2})^\intercal [{\bf u}-{\bf u}^k] \leq 0} \nonumber\\
& \added{y_1({\bf u}) := u_1^2 + u_2^2}\nonumber\\
& \added{y_2({\bf u}) := 1 - u_1 + u_2^2} \nonumber
\end{align}
\added{and the following update rule is applied:}
\begin{align}
\label{eq:MA_filter_bias}
\added{\varepsilon^{k+1}_{\delta y_i} =}\ & \added{(1-\eta) \varepsilon^{k+1}_{\delta y_i} + \eta\left[y_i^{\rm p}({\bf u}^k) - y_i({\bf u}^k)\right]} \\
\label{eq:MA_filter_grad}
\added{\boldsymbol{\lambda}_{\delta y_i}^{k+1} =}\ & \added{(1-\eta) \boldsymbol{\lambda}_{\delta y_i}^k + \eta\left[\boldsymbol{\nabla} y_i^{\rm p}({\bf u}^k) - \boldsymbol{\nabla} y_i({\bf u}^k)\right]}
\end{align}
\added{for a given gain value $0<\eta\leq 1$. For simplicity, we may determine the plant gradients $\boldsymbol{\nabla} y_i^{\rm p}({\bf u}^k)$ using forward finite differences, e.g. with steps $\Delta u_1=\Delta u_2=0.1$. The comparison of multiple modifier-adaptation runs on Figures~\ref{fig:ex1_std_08}--\ref{fig:ex1_std_02} for gain values of $\eta=0.8$, $0.5$ and $0.2$, and the corresponding cost envelopes on Figure~\ref{fig:ex1_std_obj}, confirms that this basic scheme can steer the iterates to a neighborhood of the plant optimum. The use of large gain values (e.g., $\eta=0.8$) enables a fast adaptation, but the iterates exhibit a high variance around the plant optimum. Decreasing the gain value (e.g., $\eta=0.2$) reduces this variance, yet at the cost of a significantly slower adaptation.}

\replaced{By contrast, t}{T}he modified optimization problem that is solved at each iteration \added{of Algorithm~\ref{alg:MA-GP-TR}} to determine the next \deleted{RTO} move is \deleted{thus} given by:
\begin{align}
{\bf d}^{k+1} \in \argmin_{\|{\bf d}\|\leq\Delta^{k}}~ & \mathcal{A}[y_1+\mu^k_{\delta y_1},\sigma^k_{\delta y_1},\cdot]({\bf u}^k+{\bf d}) \label{eq:Illu1_ma}\\
\text{s.t.}\ \ & [y_2+\mu^k_{\delta y_2}]({\bf u}^k+{\bf d}) \leq 0 \nonumber\\
& y_1({\bf u}) := u_1^2 + u_2^2\nonumber\\
& y_2({\bf u}) := 1 - u_1 + u_2^2 \nonumber
\end{align}
where the GP modifiers capture the output mismatch, $y_i^{\rm p}-y_i \sim \mathcal{GP}\left(\mu_{\delta y_i},\sigma_{\delta y_i}^{2}\right)$, for $i=1,2$; and the acquisition function $\mathcal{A}$ may either be LCB or EI (cf. Section~\ref{sec:GP}) or the cost predictor itself ($y_1+\mu^k_{\delta y_1}$) if exploration is not considered. \added{We analyze the effect of various designs choices in the following subsections.}

\begin{figure}[htbp]
\begin{subfigure}{0.5\textwidth}
  \centering
  \includegraphics[width=.98\linewidth]{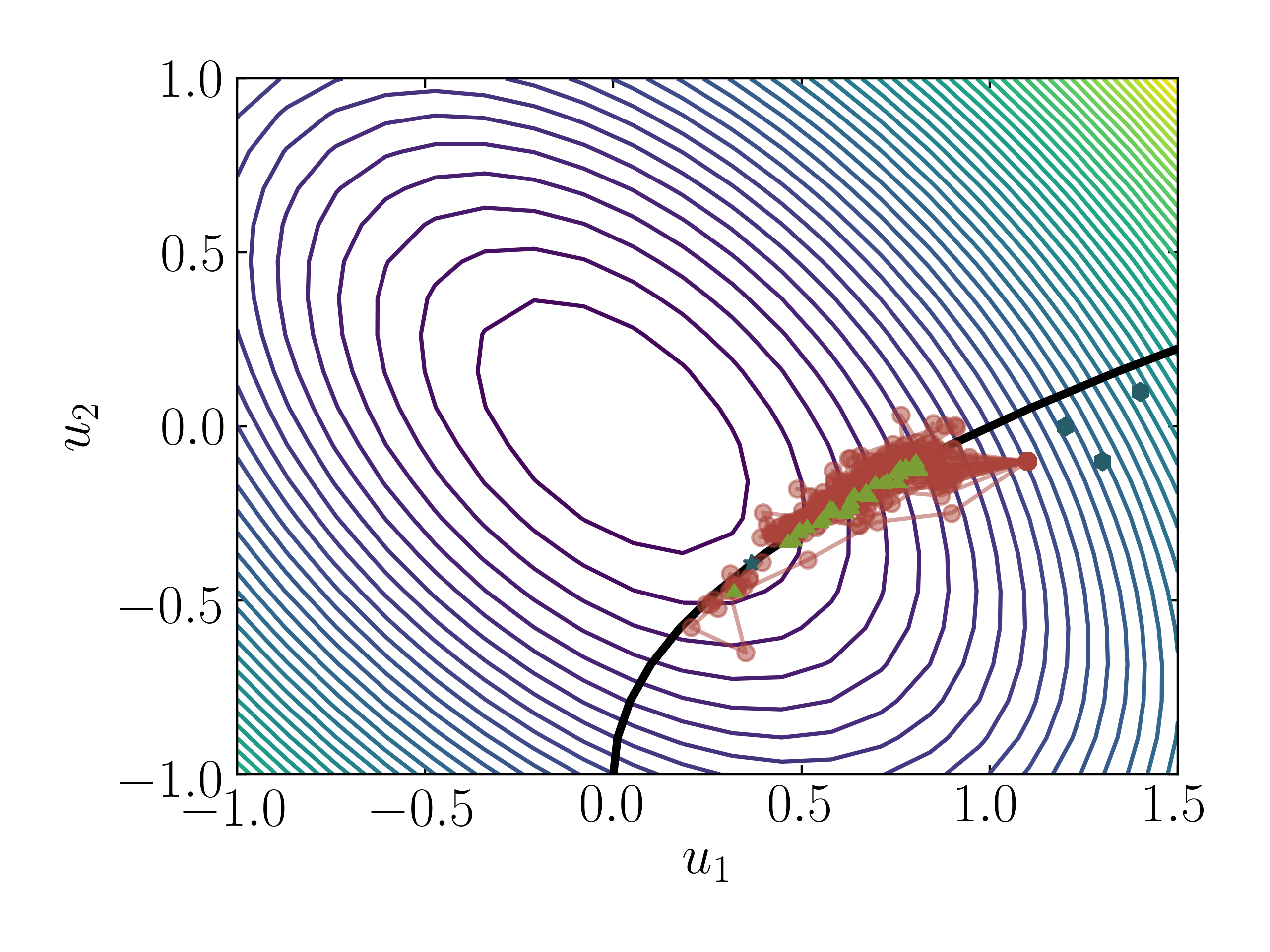}\vspace{-.85em}
  \caption{RTO iterations without acquisition function}
  \label{fig:ex1_expl_no}
\end{subfigure}
\begin{subfigure}{0.5\textwidth}
  \centering
  \includegraphics[width=.98\linewidth]{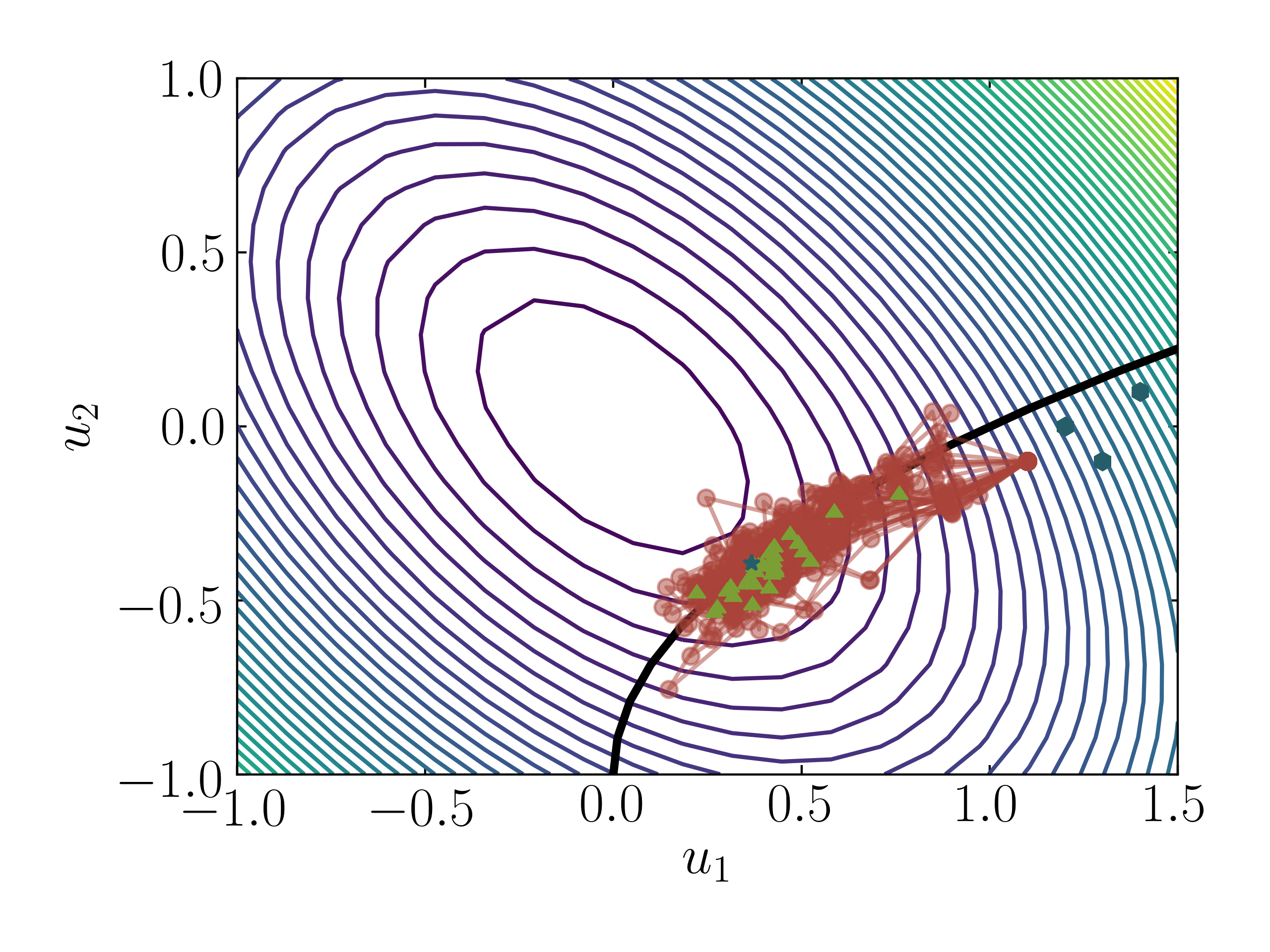}\vspace{-.85em}
  \caption{RTO iterations with LCB acquisition function}
  \label{fig:ex1_expl_ucb}
\end{subfigure}
\begin{subfigure}{0.5\textwidth}
  \centering
  \includegraphics[width=.98\linewidth]{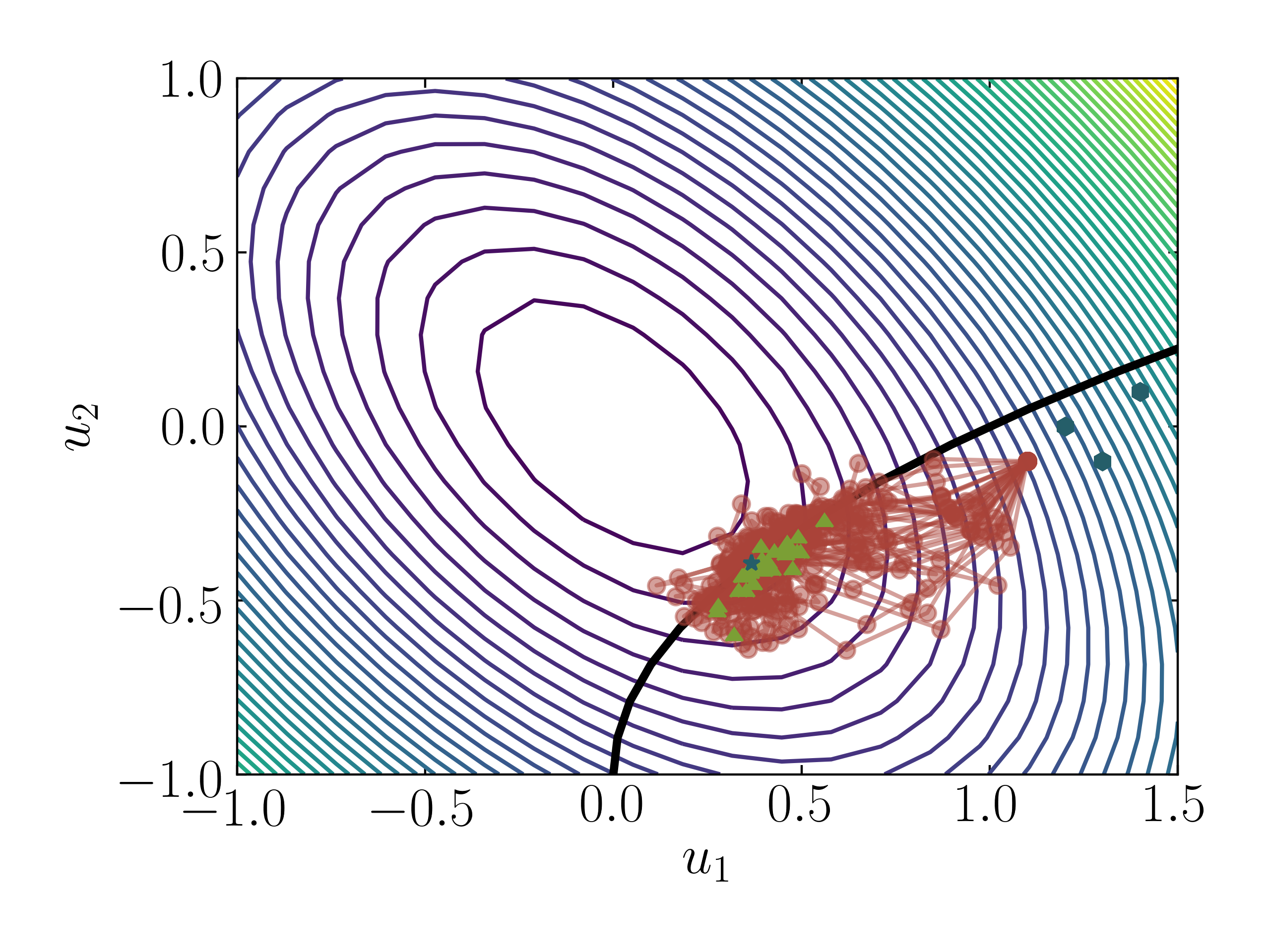}\vspace{-.85em}
  \caption{RTO iterations with EI acquisition function}
  \label{fig:ex1_expl_ei}
\end{subfigure}
\begin{subfigure}{0.5\textwidth}
  \centering
  \includegraphics[width=.98\linewidth]{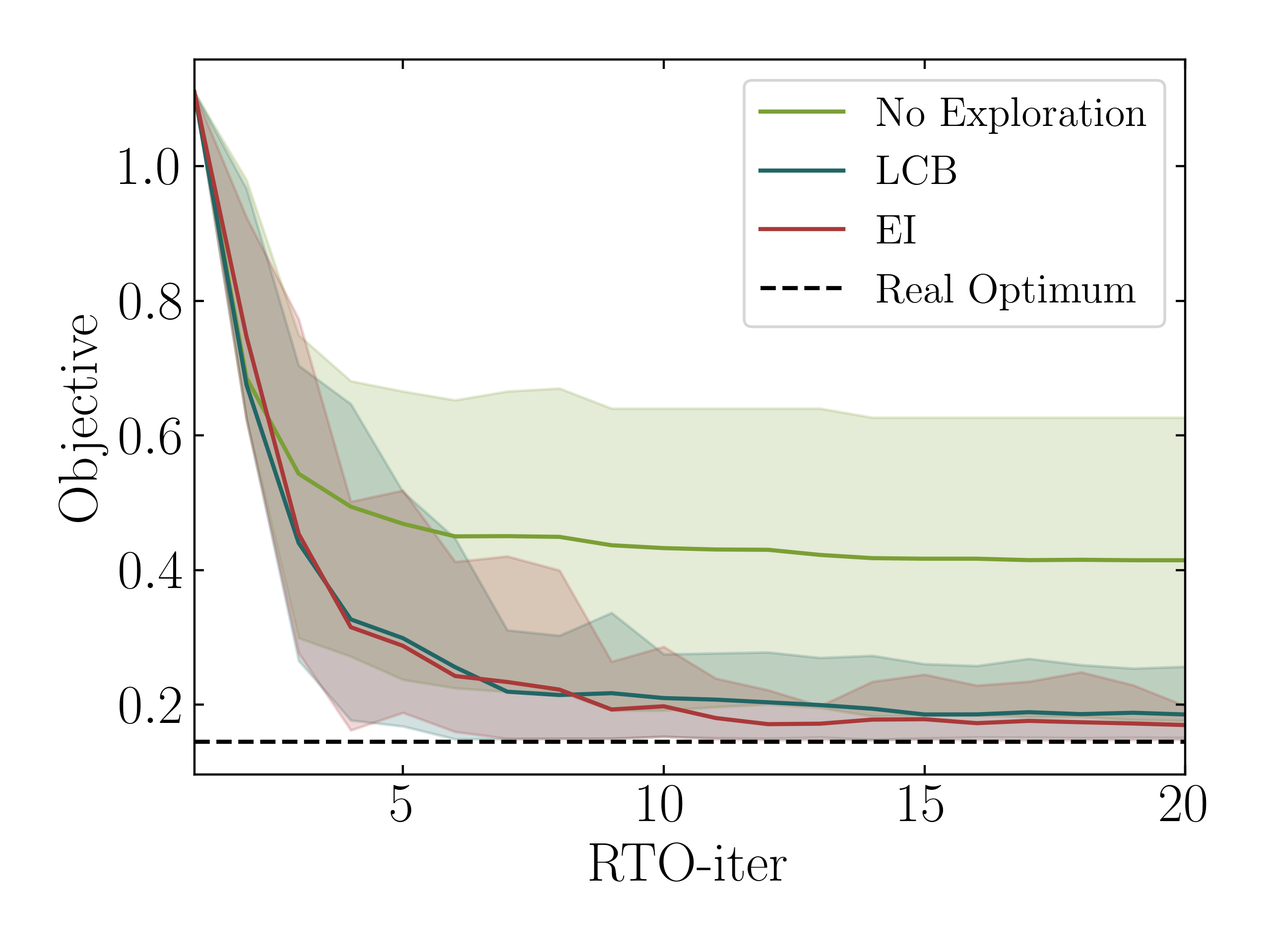}\vspace{-.85em}
  \caption{Evolution of process cost without and with acquisition function}
  \label{fig:ex1_expl_obj}
\end{subfigure}
\caption{RTO iterations for Problem~\eqref{eq:Illu1} corresponding to various exploration strategies in Algorithm~\ref{alg:MA-GP-TR}. A prior process model is used and \replaced{the process noise is estimated by the GPs}{no prior knowledge of the process noise is assumed}. {\bf (a), (b), (c)} Clouds of iterates (red connected circles) for 30 process noise realizations, initialized from the same sample points (blue hexagons) and interrupted after 20 RTO iterations (green triangles); the process optimum is depicted with a blue star. {\bf (d)}~Evolution of the 95th percentile of process cost values over all the noise realizations with the RTO iterations (showing only the feasible iterates). 
\label{fig:ex1_exploration}}
\end{figure}

\paragraph{On the Benefit of Using an Acquisition Function}

A key feature of Algorithm~\ref{alg:MA-GP-TR} lies in the use of an acquisition function for promoting exploration within the active trust region, rather than enforcing \added{(probabilistic)} full linearity of the surrogate models in Problem~\eqref{eq:modified_problem_GP+TR}. The behavior of several modifier-adaptation schemes, without and with such an acquisition function, is compared in Figure~\ref{fig:ex1_exploration} for multiple realizations of the process noise. 

The comparison of multiple modifier-adaptation runs on Figure~\ref{fig:ex1_expl_no} and the corresponding cost envelope on Figure~\ref{fig:ex1_expl_obj} clearly show that, without adding an exploration term in the modified cost of Problem~\eqref{eq:Illu1_ma}, certain RTO runs may get trapped at a suboptimal point. This behavior was not observed under noiseless conditions and is thus attributed to the presence of process noise. A possible cause could be the lack of a model-improvement step and enforcement of full linearity in Algorithm~\ref{alg:MA-GP-TR} (cf. Appendix~\ref{app:convergence}).

By contrast, with the modifier-adaptation schemes that use either the LCB or EI acquisition function (Figures~\ref{fig:ex1_expl_ucb} \& \ref{fig:ex1_expl_ei}), the iterates are much more likely to converge to the plant optimum in the presence of noise. This confirms the benefit of adding excitation in the modified cost of the RTO subproblems and that the selected acquisition functions are indeed suitable. Notice that the paths followed by the iterates of both schemes are comparable, although EI seems to drive the iterates more into the interior of the feasible region on this particular example. The comparison between LCB and EI on Figure~\ref{fig:ex1_expl_obj} also suggests that the latter may promote a faster progress and a lower variance around the plant optimum\added{, though no attempt was made to tailor the parameter $\beta$ in the LCB function \eqref{eq:LCB}.}

\added{It is also worth mentioning that the spread of the iterate clouds on Figure~\ref{fig:ex1_exploration} -- in particular at the final RTO iteration (green triangles) -- are significantly reduced compared to the basic modifier-adaptation scheme results on Figure~\ref{fig:ex1_standard}. This is attributed to the ability of the GP surrogates to detect and/or filter out noise during the regression process. This noise reduction capability is also shared with other surrogate modeling approaches in modifier-adaption schemes, such as quadratic approximations \citep{Gao2016}.}

\begin{figure}[htbp]
\begin{subfigure}{0.5\textwidth}
  \centering
  \includegraphics[width=.98\linewidth]{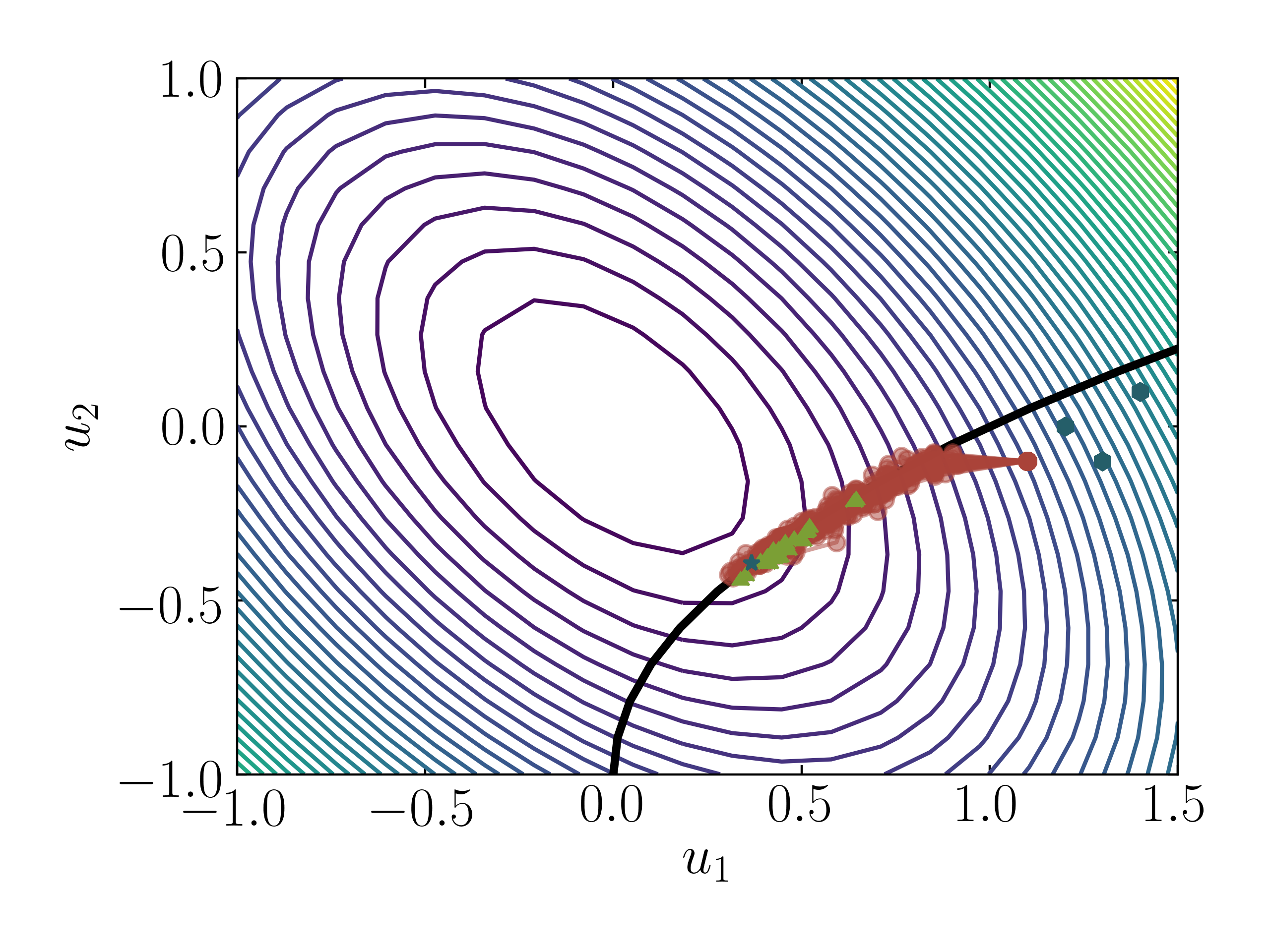}\vspace{-.85em}
  \caption{RTO iterations without acquisition function}
  \label{fig:ex1_expl2_no}
\end{subfigure}
\begin{subfigure}{0.5\textwidth}
  \centering
  \includegraphics[width=.98\linewidth]{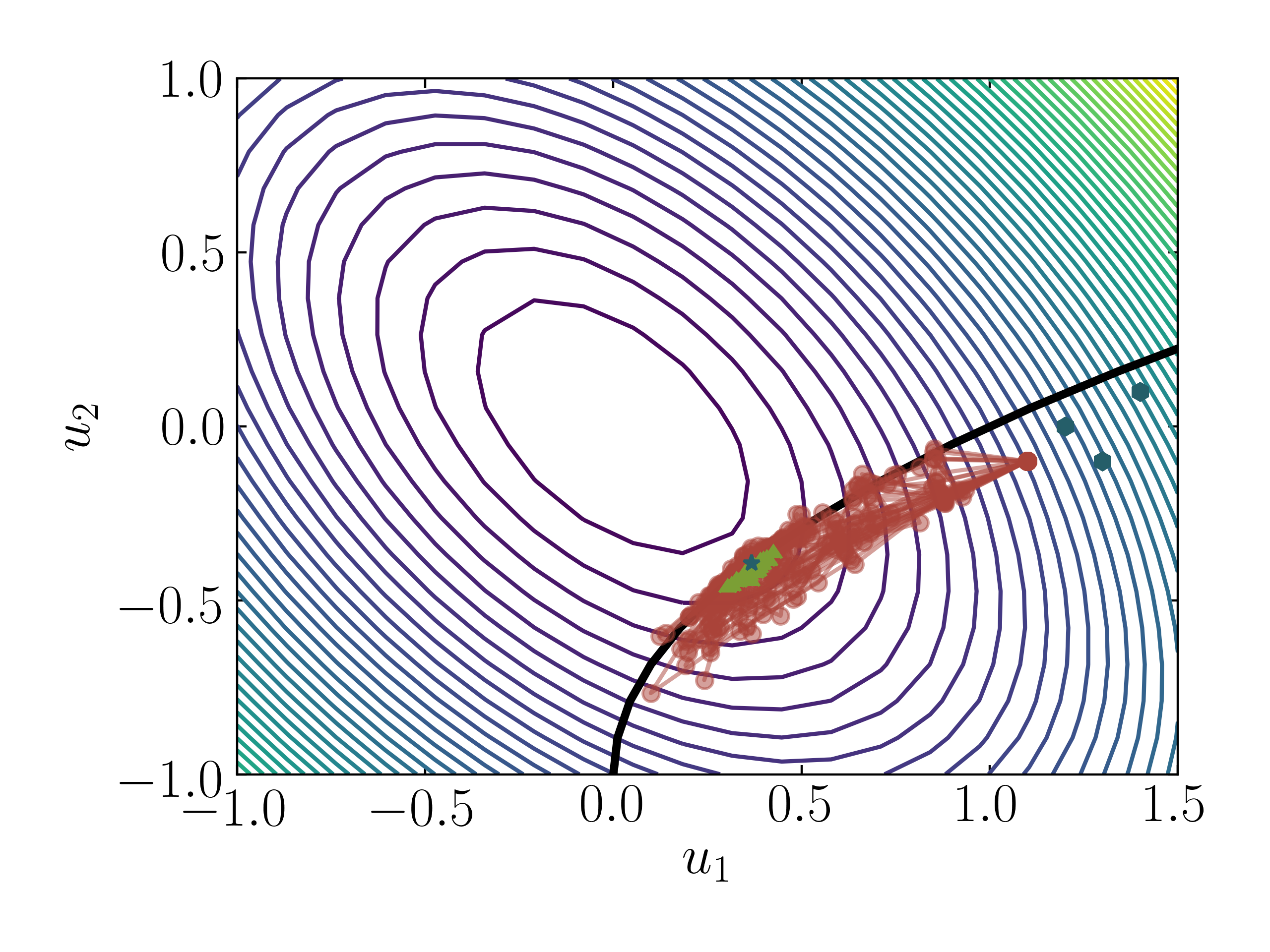}\vspace{-.85em}
  \caption{RTO iterations with LCB acquisition function}
  \label{fig:ex1_expl2_ucb}
\end{subfigure}
\begin{subfigure}{0.5\textwidth}
  \centering
  \includegraphics[width=.98\linewidth]{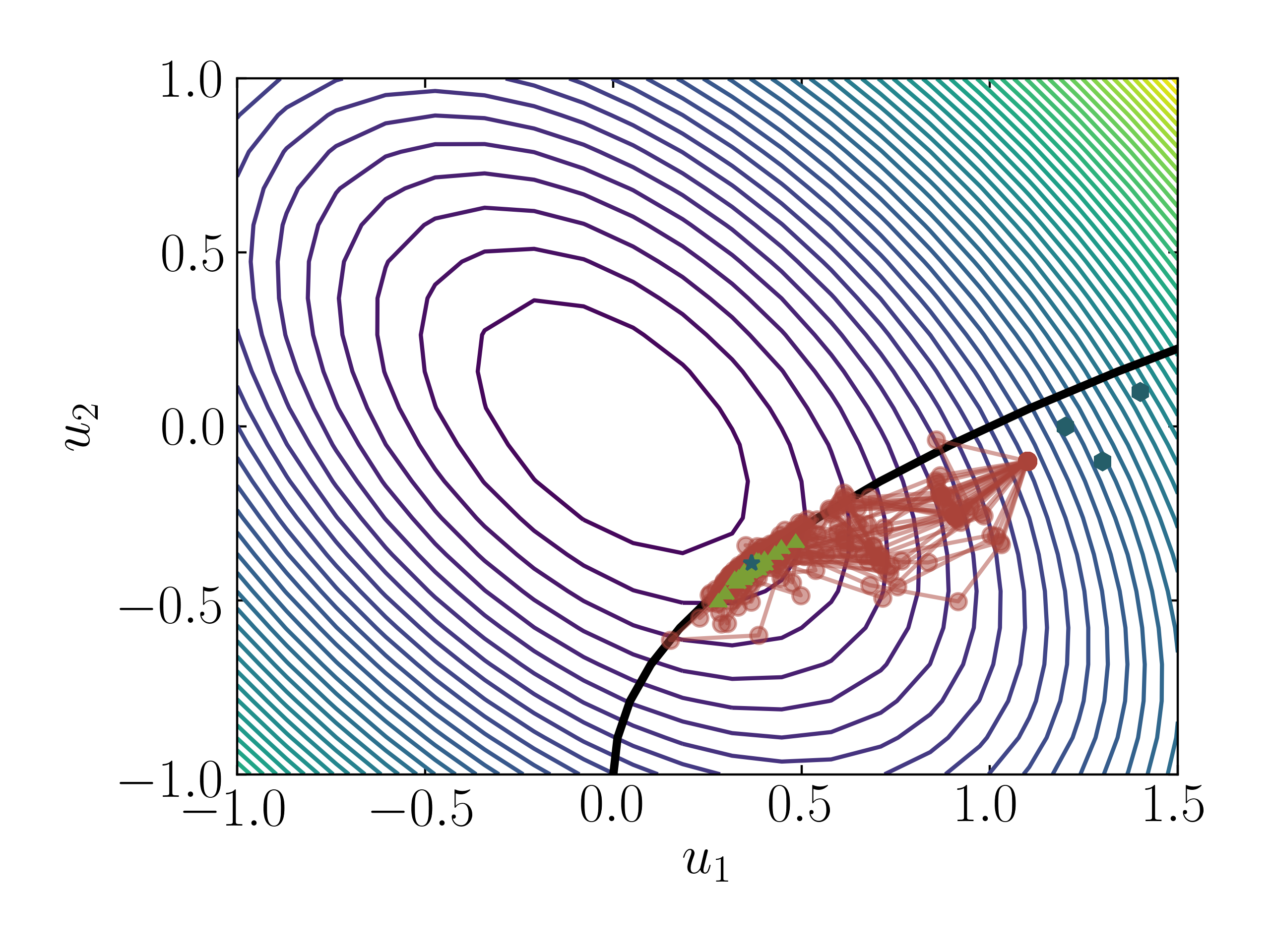}\vspace{-.85em}
  \caption{RTO iterations with EI acquisition function}
  \label{fig:ex1_expl2_ei}
\end{subfigure}
\begin{subfigure}{0.5\textwidth}
  \centering
  \includegraphics[width=.98\linewidth]{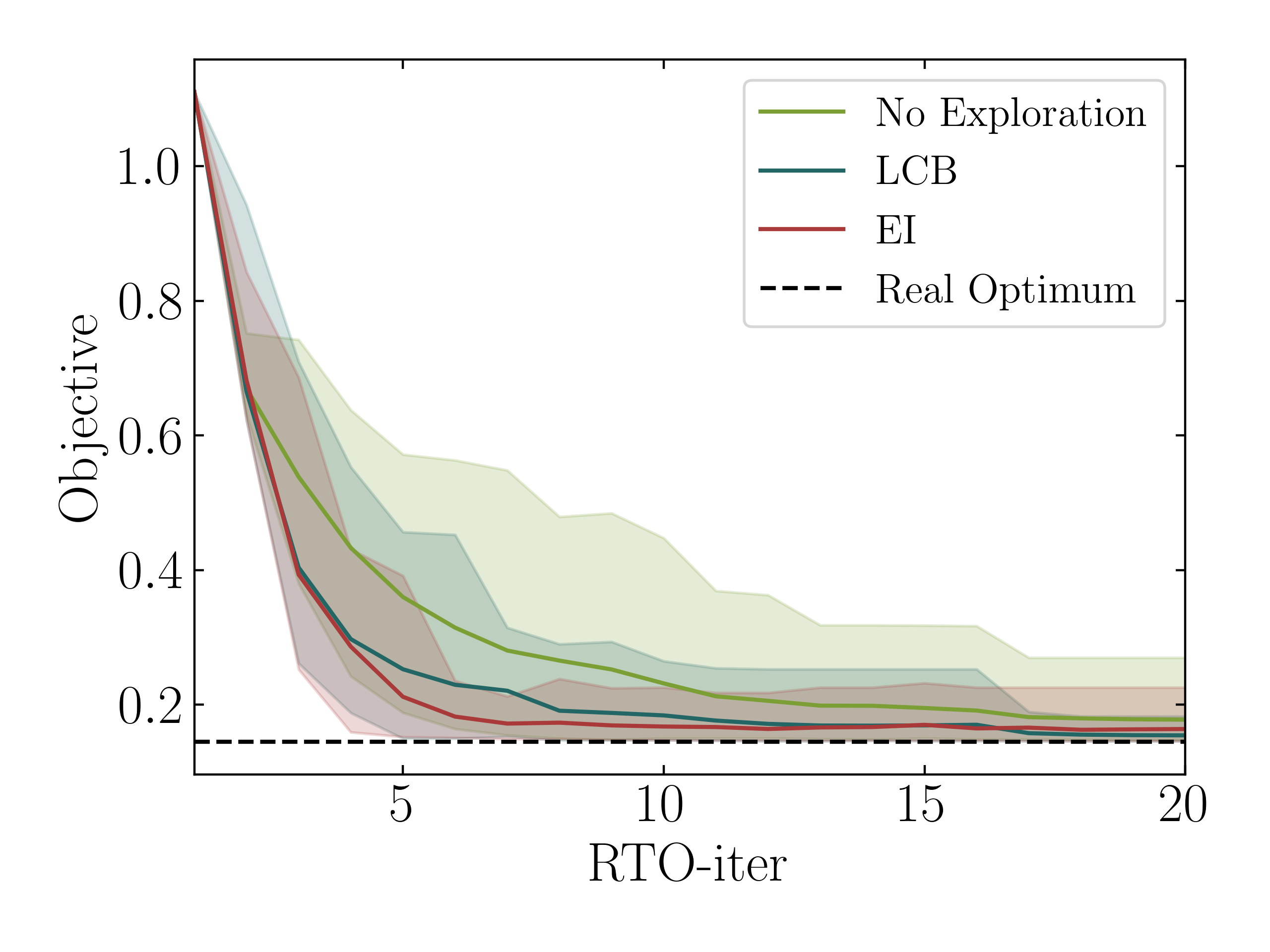}\vspace{-.85em}
  \caption{Evolution of process cost \replaced{given the correct process noise level}{without and with acquisition function}}
  \label{fig:ex1_expl2_obj}
\end{subfigure}
\begin{subfigure}{0.5\textwidth}
  \centering
  \includegraphics[width=.98\linewidth]{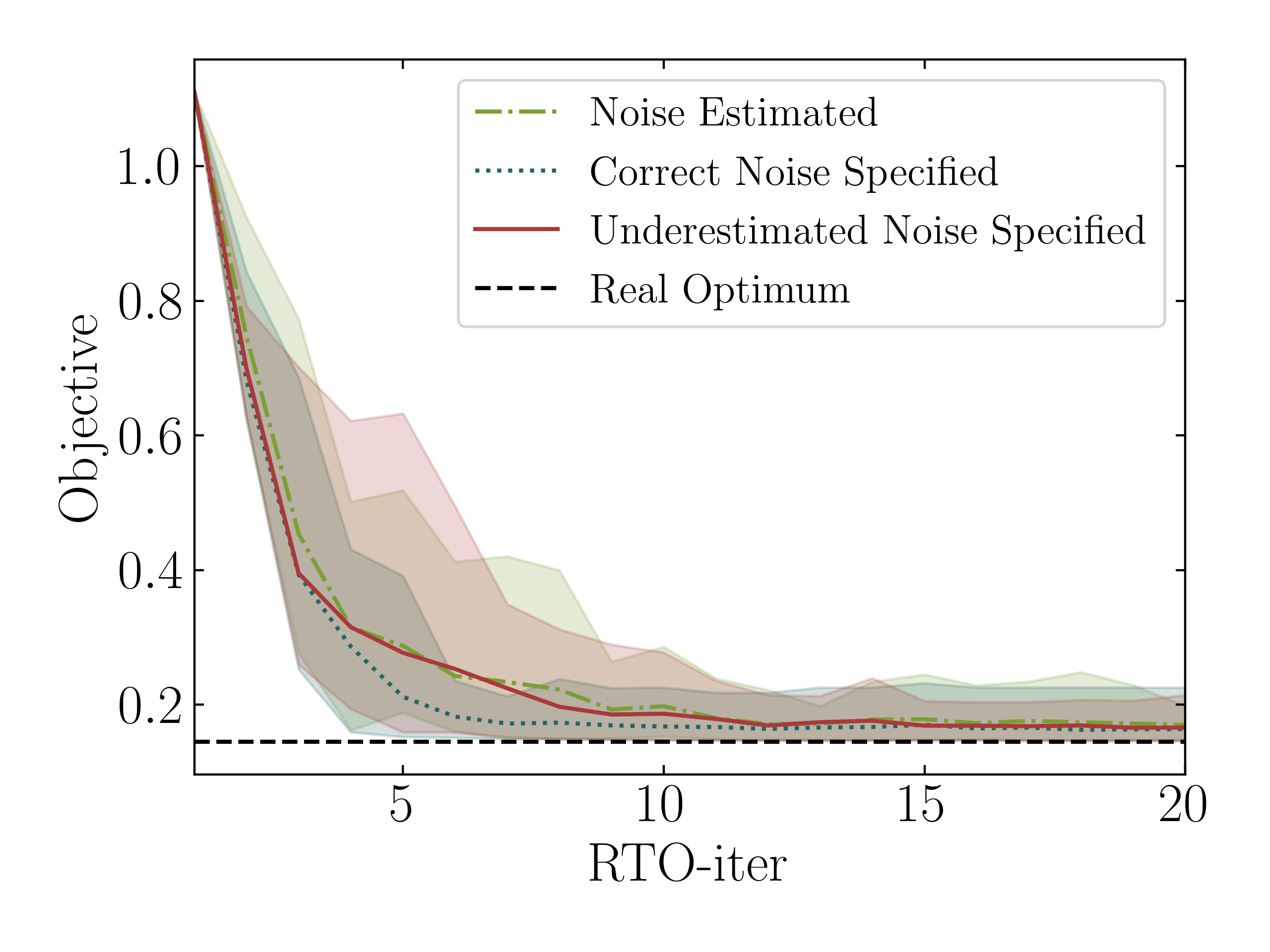}\vspace{-.85em}
  \caption{\added{Evolution of process cost with EI under various noise assumptions}}
  \label{fig:ex1_expl2_obj2}
\end{subfigure}
\caption{RTO iterations for Problem~\eqref{eq:Illu1} corresponding to various exploration \added{and noise handling }strategies in Algorithm~\ref{alg:MA-GP-TR}. A prior process model is used\deleted{ and the process noise level is assumed to be given}. {\bf (a), (b), (c)} Clouds of iterates (red connected circles) for 30 process noise realizations, initialized from the same sample points (blue hexagons) and interrupted after 20 RTO iterations (green triangles)\added{, with correct noise level specified to the GPs}; the process optimum is depicted with a blue star. {\bf (d), (e)} Evolution of the 95th percentile of process cost values over all the noise realizations with the RTO iterations (showing only the feasible iterates)\added{, corresponding to different assumptions about the noise level}. 
\label{fig:ex1_exploration2}}
\end{figure}

\paragraph{On the Benefit of \replaced{Specifying}{Knowing} the Process Noise\added{ Level}}

Not only does GP regression provide a natural approach to describing the plant-model mismatch in a non-parametric way, but it also enables estimating the variance of the observations alongside the other hyperparameters of the GPs in case the noise level is unspecified. In our initial case study shown in Figure~\ref{fig:ex1_exploration}, \replaced{the process noise was estimated in this way}{we did not assumed any prior knowledge of the process noise}. By contrast, the results in Figure~\ref{fig:ex1_exploration2} \replaced{make some prior assumptions regarding the}{assume that the (correct)} process noise \replaced{level}{is given}. 

By and large, the performance of all three modifier-adaptation schemes is clearly enhanced by the specification of the \added{correct }noise level \added{(Figure~\ref{fig:ex1_expl2_obj})}. In the case that no extra excitation is added to the modified cost of the RTO subproblems (Figure~\ref{fig:ex1_expl2_no}), the odds of the iterates getting trapped at a suboptimal point are significantly reduced, albeit still not negligible; while with an acquisition function (Figures~\ref{fig:ex1_expl2_ucb} \& \ref{fig:ex1_expl2_ei}), the variability of the iterates around the plant optimum is much lower. \added{A final comparison is conducted where the process noise level specified to the GPs is incorrect (Figure~\ref{fig:ex1_expl2_obj2}), in this case underestimating the noise variance by a factor of $2$. Such a misspecification negates the benefits of providing the noise level, showing a comparable performance as with noise estimation during GP training. A larger underestimation of the noise level could even become detrimental to the RTO system's reliability.} This cursory \replaced{analysis}{comparison} illustrates well the benefits of characterizing the process noise, e.g. based on historical data.

\paragraph{On the Benefit of Specifying a Nominal Process Model}

\begin{figure}[htb]
\begin{subfigure}{0.5\textwidth}
  \centering
  \includegraphics[width=.98\linewidth]{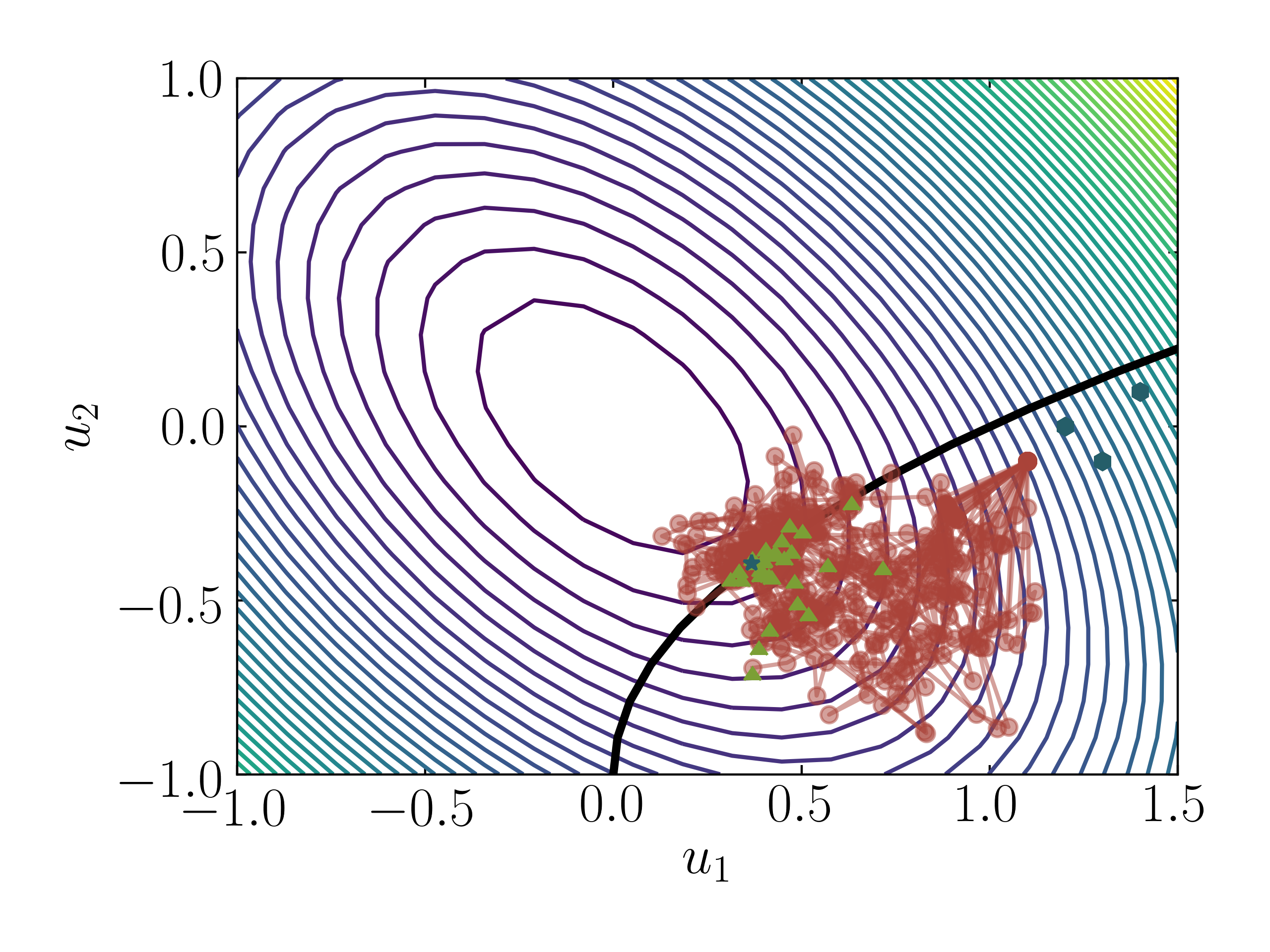}\vspace{-.85em}
  \caption{RTO iterations with EI acquisition function}
  \label{fig:ex1_expl3_ei}
\end{subfigure}
\begin{subfigure}{0.5\textwidth}
  \centering
  \includegraphics[width=.98\linewidth]{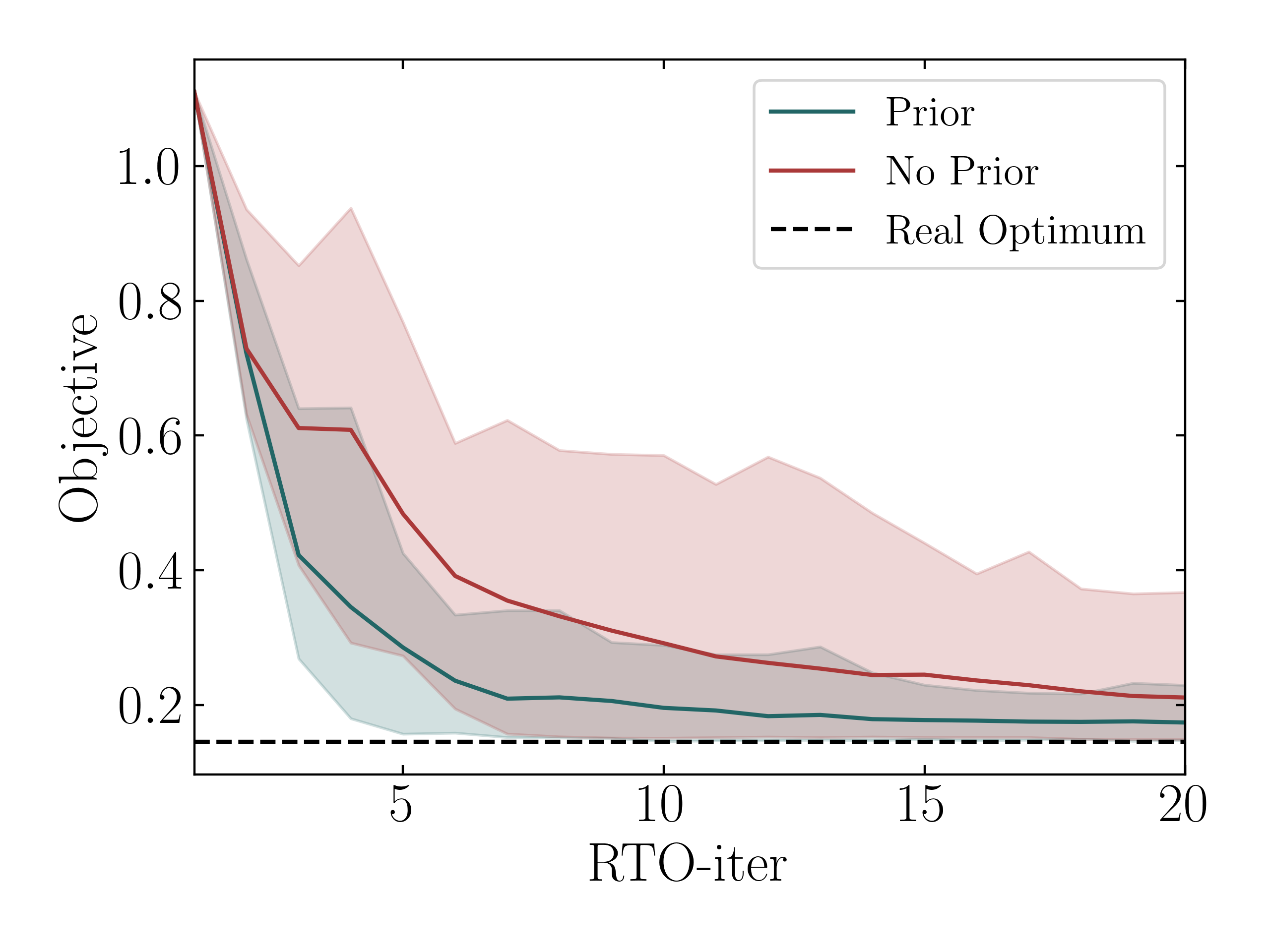}\vspace{-.85em}
  \caption{Evolution of process cost without and with a prior process model}
  \label{fig:ex1_expl3_obj}
\end{subfigure}
\caption{RTO iterations for Problem~\eqref{eq:Illu1} generated by Algorithm~\ref{alg:MA-GP-TR} without and with a prior model. The EI acquisition function is used and no prior knowledge of the process noise is assumed. {\bf (a)} Clouds of iterates (red connected circles) for 30 process noise realizations, initialized from the same sample points (blue hexagons) and interrupted after 20 RTO iterations (green triangles); the process optimum is depicted with a blue star. {\bf (b)} 95th percentile of the process cost values over all noise realizations at each RTO iteration (showing only the feasible iterates). 
\label{fig:ex1_exploration3}}
\end{figure}

The basic idea behind modifier adaptation entails correcting a model-based optimization problem so that its solution will match the plant optimum upon convergence. By contrast, Bayesian optimization and derivative-free optimization do not rely on a preexisting model, so it seems legitimate to raise the question whether a GP model alone would be suitable to drive such an RTO system. Discarding the process model altogether is akin to model-free RTO, which is be easier to design and maintain, but could result in large performance loss or lesser reliability compared to model-based RTO nonetheless. The behavior of a modifier-adaptation scheme without a prior (nominal) model---that is, setting $y_1({\bf u})=y_2({\bf u})=0$ in Problem~\eqref{eq:Illu1_ma}---is shown in Figure~\ref{fig:ex1_exploration3}.

The cloud of iterates on Figure~\ref{fig:ex1_expl3_ei} presents a much \replaced{wider spread}{larger span} than its counterpart on Figure~\ref{fig:ex1_expl_ei} which uses a prior model and the same EI acquisition function. This behavior is also observed on Figure~\ref{fig:ex1_expl3_obj} where the envelope of cost values for a range of noise scenarios is two to three times wider after discarding the nominal process model. Many more infeasible iterates are furthermore generated in this latter scenario, which requires backing-off more frequently and thereby slows down the adaptation. The fact that several final iterates are not on the constraint in Figure~\ref{fig:ex1_expl3_ei} suggests that, without building on a prior model, the GP surrogates yield an inaccurate prediction of the actual process constraint. The reason for this could be the lack of exploration of the feasible region, since an acquisition function is only used to promote exploration in the objective function of Problem~\eqref{eq:Illu1_ma}. Improved modifier-adaptation schemes that add excitation to both the cost and constraint functions will be investigated in future work.

Overall, the comparisons conducted in this section have provided compelling evidence that (i) using an acquisition function, (ii) knowing the process noise level, and (iii) specifying a nominal process model can greatly enhance the reliability of a modifier-adaptation scheme based on GP modifiers. Naturally, the extent to which such design choices will improve an RTO system is largely problem dependent. The following section presents further results for two numerical case studies.

\section{Case Studies}\label{sec:case}

\subsection{Williams-Otto Benchmark Problem}\label{sec:case_two}

We first assess the proposed modifier-adaptation algorithm with Gaussian process, trust region and acquisition function (Algorithm~\ref{alg:MA-GP-TR}) on the classical Williams-Otto benchmark problem. A continuous stirred-tank reactor (CSTR) is fed with two streams of pure components {\sf A} and {\sf B}, with respective mass flowrates $F_{\sf A}$ and $F_{\sf B}$. The reactor operates at steady state and under the temperature $T_{\rm r}$. The chemical reactions between these reagents produce two main products {\sf P} and {\sf E}, through a series of chemical reactions that also produce an intermediate {\sf C} and a byproduct {\sf G}:
\begin{align*} 
{\sf A} + {\sf B} \longrightarrow~& {\sf C}\\
{\sf B} + {\sf C} \longrightarrow~& {\sf P} + {\sf E}\\
{\sf C} + {\sf P} \longrightarrow~& {\sf G}
\end{align*}
Structural plant-model mismatch is introduced in the problem by assuming that the approximate kinetic model only knows about the following two reactions, which omit the intermediate species~{\sf C}:
\begin{align*} 
{\sf A} + 2{\sf B} \longrightarrow~& {\sf P} + {\sf E}\\
{\sf A} + {\sf B} + {\sf P} \longrightarrow~& {\sf G}
\end{align*}
The complete set of mass-balance equations and kinetic rate equations for both reaction systems are the same as those reported by \citet{Mendoza2016} and not reproduced here for brevity.

\begin{figure}[tbp]
\begin{subfigure}{0.5\textwidth}
  \centering
  \includegraphics[width=.98\linewidth]{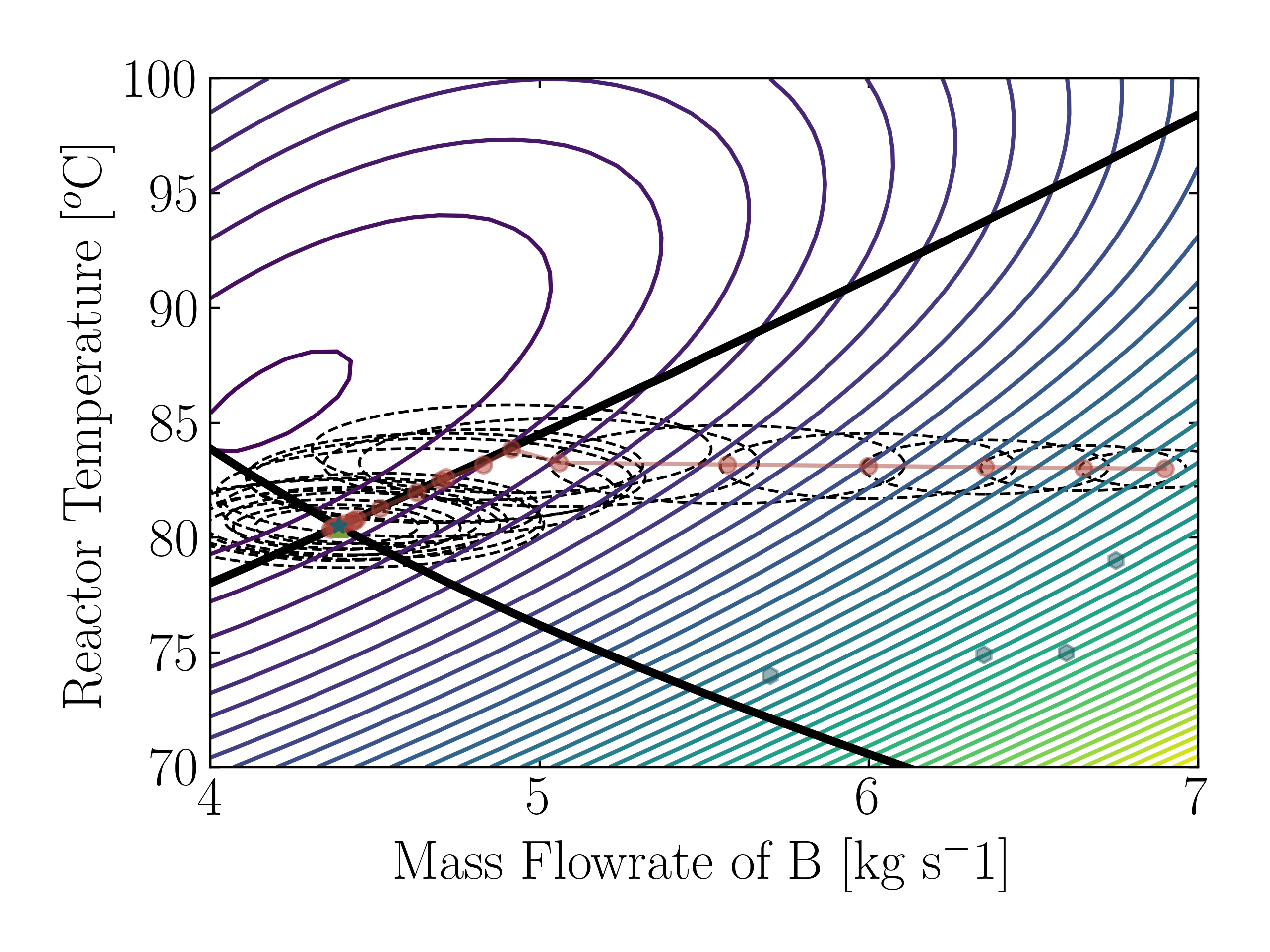}\vspace{-.85em}
  \caption{Trust regions evolution along a single RTO path}
  \label{fig:WO_pathTR}
\end{subfigure}
\begin{subfigure}{0.5\textwidth}
  \centering
  \includegraphics[width=.98\linewidth]{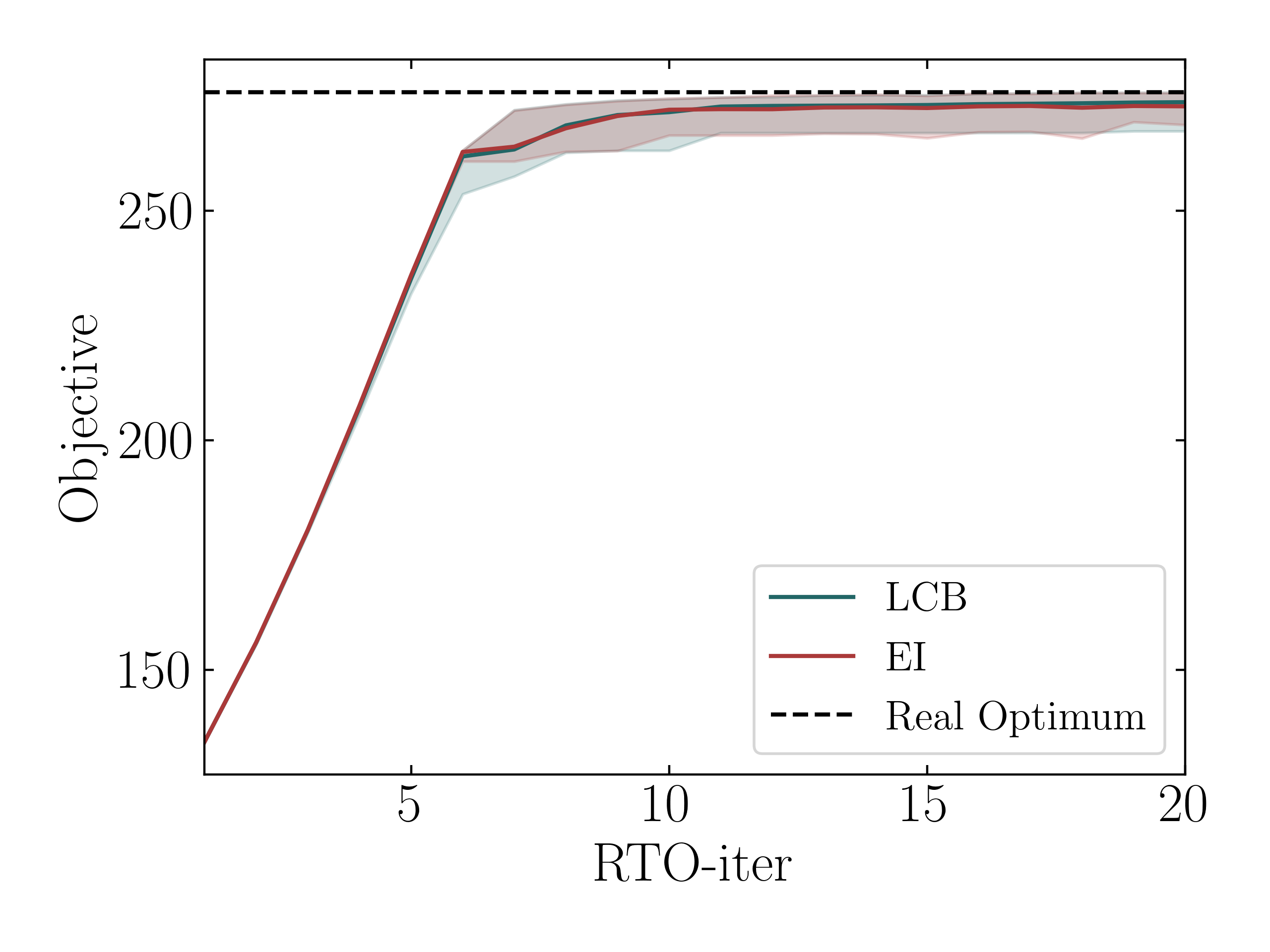}\vspace{-.85em}
  \caption{Evolution of process cost for multiple runs}
  \label{fig:WO_cost}
\end{subfigure}
\begin{subfigure}{0.5\textwidth}
  \centering
  \includegraphics[width=.98\linewidth]{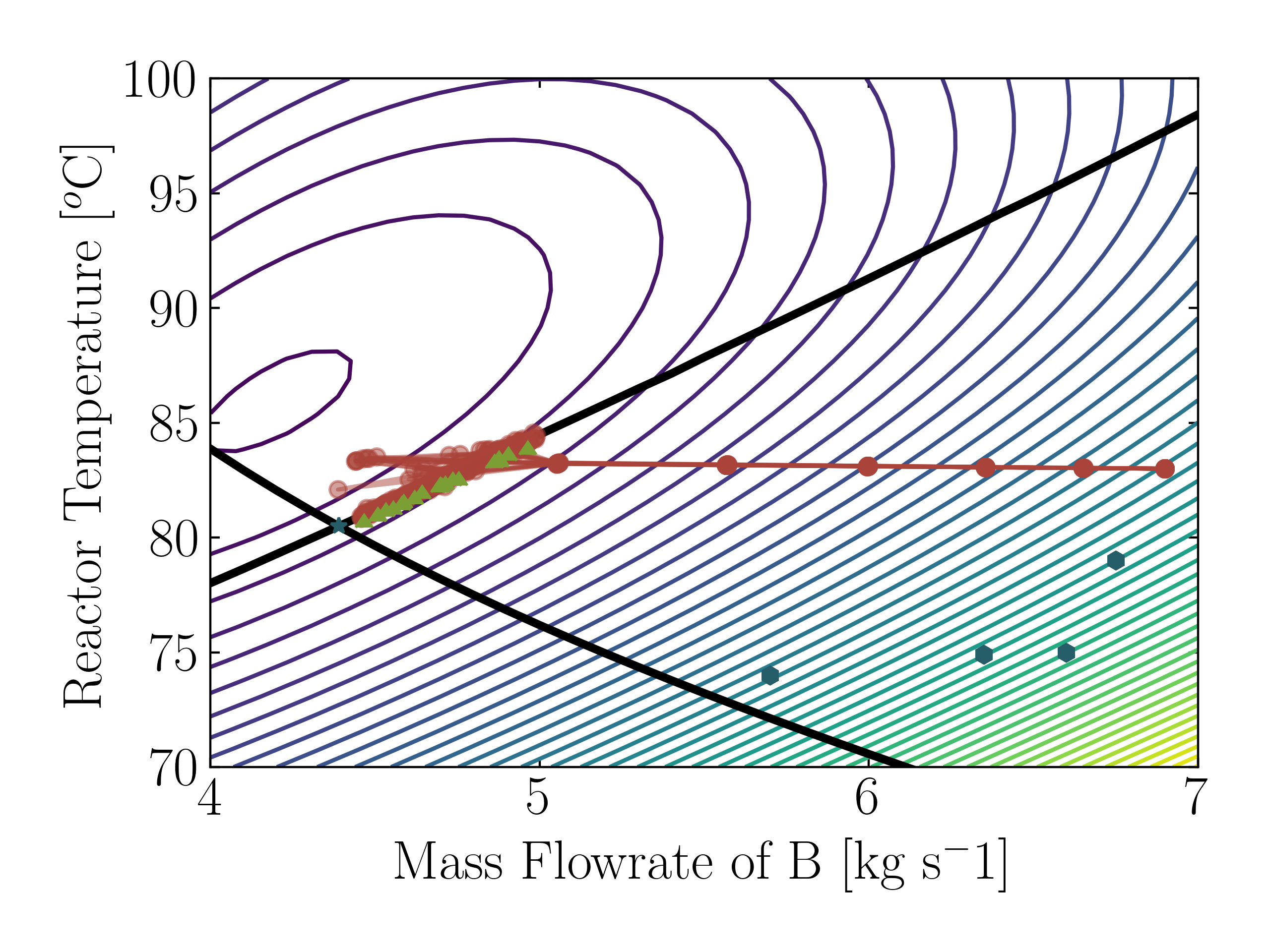}\vspace{-.85em}
  \caption{RTO iterates with EI acquisition function}
  \label{fig:WO_EI}
\end{subfigure}
\begin{subfigure}{0.5\textwidth}
  \centering
  \includegraphics[width=.98\linewidth]{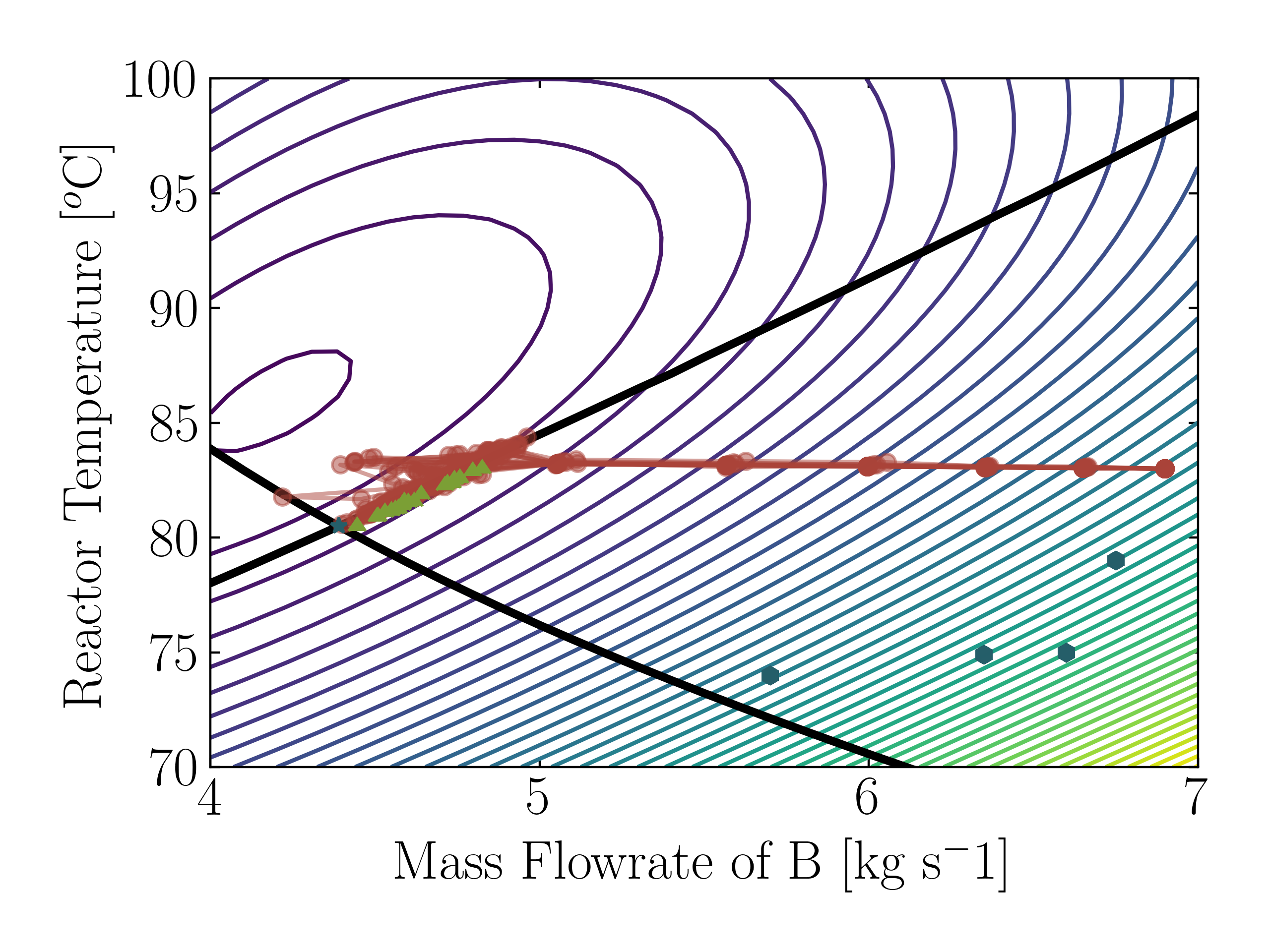}\vspace{-.85em}
  \caption{RTO iterates with LCB acquisition function}
  \label{fig:WO_LCB}
\end{subfigure}
\caption{RTO results for the Williams-Otto case study (Problem~\ref{eq:WO}) using Algorithm~\ref{alg:MA-GP-TR}. {\bf (a)} Evolution of the trust-region size (dashed ellipsoids) for a single RTO run with EI acquisition function, interrupted after 20 iterations. {\bf (b)} Evolution of the 95th percentile of process cost values over 30 noise realizations with the RTO iterations (showing only the feasible iterates). {\bf (c), (d)} Clouds of iterates (red connected circles) for 30 process noise realizations, initialized from the same sample points (blue hexagons) and interrupted after 20 RTO iterations (green triangles); the process optimum is depicted with a blue star.\label{fig:WO}}
\end{figure}

The optimization problem seeks to maximize the economic profit by manipulating the feedrate $F_{\sf B}$ and the reactor temperature $T_{\rm r}$, subject to operating constraints on the residual mass fractions of {\sf A} and {\sf G} at the reactor outlet:
\begin{align}
\min_{F_{\sf B},T_{\rm r}}~~ & G_0 := (1043.38 X_{\sf P} + 20.92 X_{\sf E})\:(F_{\sf A}+F_{\sf B}) - 79.23 F_{\sf A} - 118.34 F_{\sf B}\label{eq:WO}\\
\text{s.t.}~~ & \text{CSTR model \citep{Mendoza2016}}\nonumber\\
& G_{1} := X_{\sf A}-0.12 \leq 0\nonumber\\
& G_{2} := X_{\sf G}-0.08 \leq 0\nonumber\\
& F_{\sf B}\in[4,7],~~T_{\rm r}\in[70,100]\nonumber
\end{align}
where $X_i$ denotes the mass fraction of species $i$. A graphical depiction of the problem \eqref{eq:WO} is presented in Figure~\ref{fig:WO}, where both the contour levels of the plant cost (thin multicolored lines) and the plant constraint limits (thick black lines) are shown. The case study furthermore assumes that measurements for the cost and constraint functions are available, corrupted by Gaussian distributed noise with zero mean and standard deviation $\sigma_{G_0}= 0.5$, $\sigma_{G_1}=\sigma_{G_2}= 0.0005$. However, no prior knowledge of this noise level is assumed during the construction of the GP surrogates.

\nomenclature[A]{NLP}{nonlinear programming}

The python code used to solve this case study is made available as part of the Supporting Information. The NLP solver IPOPT \citep{IPOPT}  is used to solve the optimization subproblems in the modifier-adaptation scheme. It is combined with a simple multistart heuristic (20 random starting points) to overcome numerical failures of the NLP solver and reduce the likelihood of converging to a local optimum. 

An illustration of the trust-region evolution along a particular RTO run is presented in Figure~\ref{fig:WO_pathTR}. During the first few iterations Algorithm~\ref{alg:MA-GP-TR} follows a straight path and increases the trust-region radius $\Delta$, until the boundary of the feasible domain is reached. After that, the iterates follow the active constraint and the trust-region radius is reduced to prevent constraint violations. Here, both constraints are considered unrelaxable ($\mathcal{UC}=\{1,2\}$) and $\Delta$ is reduced by a factor of 0.8 in Step~4 after back-tracking from any infeasible move. The iterates reach a close neighborhood of the plant optimum where both constraints are active after about 10 iterations. 

A comparison between multiple modifier-adaptation runs with either the LCB or EI acquisition function is presented in Figure~\ref{fig:WO_cost}. The performance is comparable and all the runs reach a neighborhood of the plant optimum within 7--11 iterations, after which they remain in the level of noise. The corresponding clouds of iterates on Figures \ref{fig:WO_EI} \& \ref{fig:WO_LCB} confirm this rapid convergence, despite several constraint violations during the search. Some of the final points after 20 iterations (green triangles) appear to be quite distant from the plant optimum, which is caused by the low sensitivity of the cost along one of the active constraints in comparison to the noise level; that is, the iterates do not get stuck at a suboptimal point.

Finally, it is worth pointing out that the performance of Algorithm~\ref{alg:MA-GP-TR} on this benchmark problem, both in terms of speed and reliability, is comparable that of other modifier-adaptation schemes. This includes the approach by \citet{Gao2016} which combines modifier adaptation with quadratic surrogates and the nested modifier-adaptation approach by \citet{Navia2015}.

\subsection{Batch-to-Batch Bioreactor Optimization}\label{sec:case_batchbio}

Our final case study investigates the performance of the proposed methodology in higher-dimensional RTO problems. We consider the batch-to-batch optimization of a photobioreactor for the production of phycocyanin ({\sf P}) by the blue-green cyanobacterium {\em Arthrospira platensis} ({\sf X}) growing on nitrates ({\sf N}). A dynamic model describing the concentrations $C_{\sf X}\ \rm [g\,L^{-1}]$, $C_{\sf N}\ \rm [mg\,L^{-1}]$ and $C_{\sf P}\ \rm [mg\,L^{-1}]$ in the  photobioreactor is given by \citep{Bradford2020}: 
\begin{align}
\dot{C}_{\sf X}(t) =\ & u_{\sf m} \dfrac{I(t)}{I(t) + k_{\sf s} + I(t)^2/k_{\sf i}} \dfrac{ C_{\sf N}(t)}{C_{\sf N}(t) + K_{\sf N}}C_{\sf X}(t) - u_{\sf d}C_{\sf X}(t)\label{ode:X}\\
\dot{C}_{\sf N}(t) =\ & -Y_{\sf N/X} u_{\sf m} \dfrac{I(t)}{I(t) + k_{\sf s} + I(t)^2/k_{\sf i}} \dfrac{ C_{\sf N}(t)}{C_{\sf N}(t) + K_{\sf N}}C_{\sf X}(t) + F_{\sf N}(t) \label{ode:N}\\
\dot{C}_{\sf P}(t) =\ & k_{\sf m} \dfrac{I(t)}{I(t) + k_{\sf sq} + I(t)^2/k_{\sf iq}}{C}_{\sf X}(t) - k_{\sf d} \dfrac{C_{\sf P}(t)}{C_{\sf N}(t) + K_{\sf Np}} \label{ode:P}
\end{align}
where the light intensity $I(t)$ $\rm [\mu E\,m^2\,s^{1}]$ and the nitrate inflow rate $F_{\sf N}(t)$ $\rm [mg\,L^{-1}\,h^{-1}]$ are manipulated inputs; and the values of the model parameters $k_{\sf d}$, $k_{\sf m}$, $k_{\sf s}$, $k_{\sf i}$, $k_{\sf sq}$, $k_{\sf iq}$, $K_{\sf N}$, $K_{\sf Np}$, $u_{\sf d}$, $u_{\sf m}$, $Y_{\sf N/X}$ are the same as those reported by \citet{Bradford2020}. For simplicity, the mass-balance equations \eqref{ode:X}--\eqref{ode:P} neglect the change in volume due to the nitrate addition and the kinetic model assumes nutrient-replete growth conditions.

The optimization problem seeks to maximize the end-batch concentration of phycocyanin after 240 hours of operation. Regarding constraints, the phycocyanin-to-cyanobacterial-biomass ratio must be kept under $1.1\ \rm wt\%$ at all times; the nitrate concentration must be kept under $800\ \rm mg\,L^{-1}$ at all times and below $150\ \rm mg\,L^{-1}$ at the end of the batch; and both manipulated inputs are bounded. A mathematical formulation of this (dynamic) optimization problem is as follows:
\begin{align}
\min_{I(t),F_{\sf N}(t)}~~ & C_{\sf P}(240)\label{eq:PBR}\\
\text{s.t.}~~ & \text{PBR model \eqref{ode:X}--\eqref{ode:P}}\nonumber\\
& C_{\sf X}(0) = 1,\ C_{\sf N}(0) = 150,\ C_{\sf P}(0) = 0\nonumber\\
& C_{\sf P}(t) \leq 0.011 C_{\sf X}(t),\ \forall t\nonumber\\
& C_{\sf N}(t) \leq 800,\ \forall t\nonumber\\
& C_{\sf N}(240) \leq 150\nonumber\\
& 120 \leq I(t) \leq 400,\ \forall t\nonumber\\
& 0 \leq F_{\sf N}(t) \leq 40,\ \forall t\nonumber
\end{align}
In order to recast it as a finite-dimensional optimization problem, both control trajectories are discretized using a piecewise-constant parameterization over 6 equidistant stages (of 60 hours each). The batch-to-batch optimization therefore comprises a total of 12 degrees of freedom. The state path constraints are also discretized and enforced at the end of each control stage.

The case study assumes that the concentrations $C_{\sf X}$, $C_{\sf N}$ and $C_{\sf P}$ can all be measured during or at the end of the batch as necessary. Process noise is simulated in this virtual reality by adding a Gaussian white noise with zero mean and standard deviation $\sigma_{C_{\sf X}}= 0.02\ \rm [g\,L^{-1}]$, $\sigma_{C_{\sf N}}= 0.316\ \rm [mg\,L^{-1}]$, and $\sigma_{C_{\sf P}}= 0.0001\ \rm [mg\,L^{-1}]$. However, no prior knowledge of this measurement noise is assumed during the construction of the GP surrogates for the cost and constraint defects. We also depart from the previous case studies by using a Mat\'ern kernel (with parameter $\nu=\frac{3}{2}$) instead of the usual squared-exponential kernel (cf. Section~\ref{sec:GP}).

Next, Algorithm~\ref{alg:MA-GP-TR} is applied to solve Problem~\eqref{eq:PBR}, both without and with the use of a nominal process model. The following dynamic model is used for the latter, which presents a structural mismatch with the plant model \eqref{ode:X}--\eqref{ode:P} regarding the light inhibition kinetics:
\begin{align}
\dot{C}_{\sf X}(t) =\ & u_{\sf m} \dfrac{I(t)}{I(t) + k_{\sf s}} \dfrac{ C_{\sf N}(t)}{C_{\sf N}(t) + K_{\sf N}}C_{\sf X}(t) - u_{\sf d}C_{\sf X}(t)\label{ode:Xmod}\\
\dot{C}_{\sf N}(t) =\ & -Y_{\sf N/X} u_{\sf m} \dfrac{I(t)}{I(t) + k_{\sf s}} \dfrac{ C_{\sf N}(t)}{C_{\sf N}(t) + K_{\sf N}}C_{\sf X}(t) + F_{\sf N}(t) \label{ode:Nmod}\\
\dot{C}_{\sf P}(t) =\ & k_{\sf m} \dfrac{I(t)}{I(t) + k_{\sf sq}} \dfrac{c_{\sf N}}{C_{\sf N}(t) + K_{\sf N}}C_{\sf X}(t) - k_{\sf d} \dfrac{C_{\sf P}(t)}{C_{\sf N}(t) + K_{\sf Np}} \label{ode:Pmod}
\end{align}
For their numerical solutions, the resulting dynamic optimization subproblems are discretized using a 4th-order Runge Kutta scheme over 25 subintervals for each control stage. All of the NLP problems are solved using IPOPT \citep{IPOPT} interfaced with CasADi \citep{Andersson2019} for computing the required derivatives. A simple multistart heuristic (20 random starting points) is applied to overcome the numerical failures of the NLP solver and reduce the likelihood of converging to a local optimum. The python code used to solve this case study is also made available as part of the Supporting Information.

\begin{figure}[tbp]
\begin{subfigure}{0.5\textwidth}
  \centering
  \includegraphics[width=.98\linewidth]{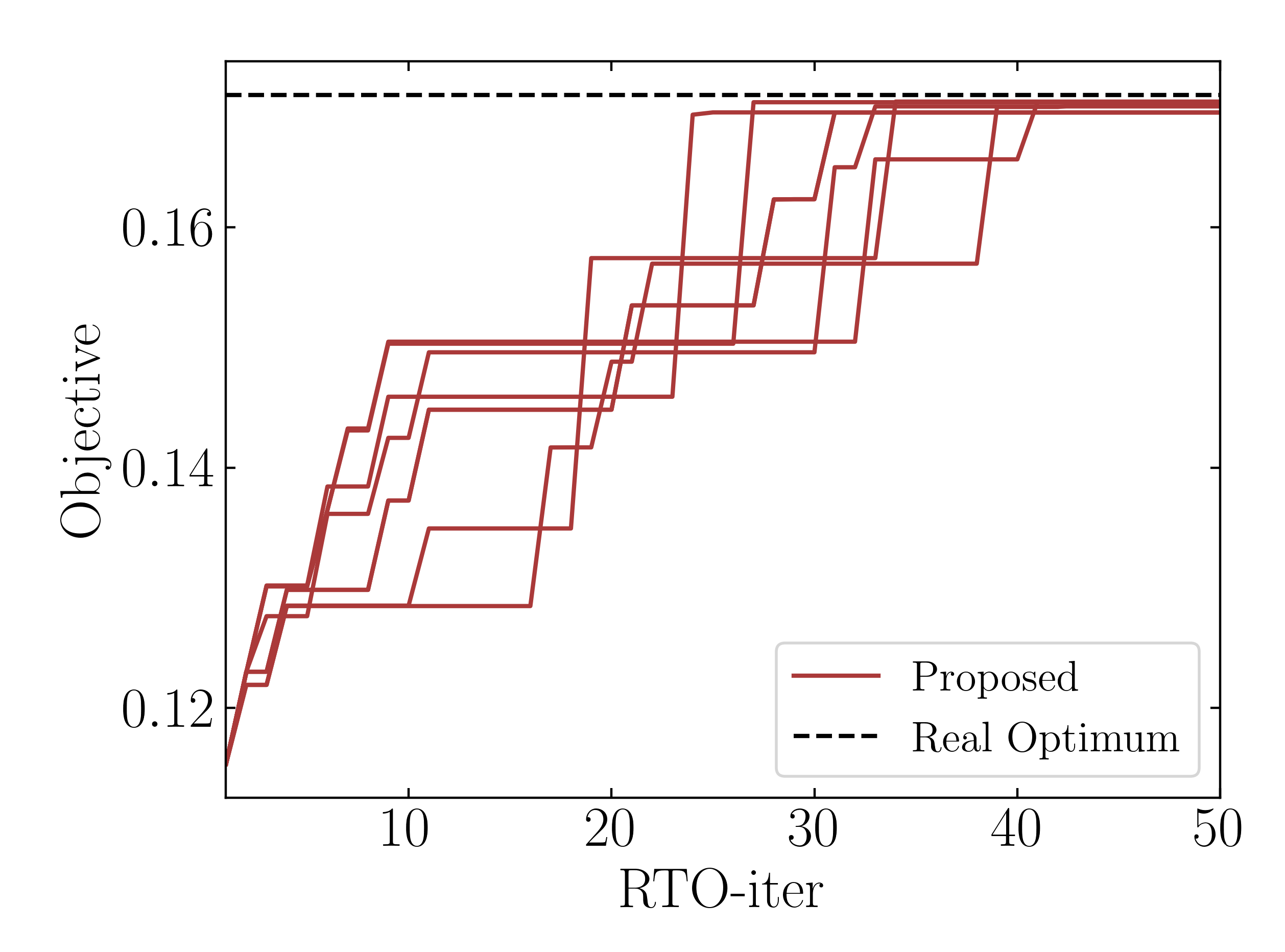}\vspace{-.5em}
  \caption{RTO iterates with prior model}
  \label{fig:PBR_prior}
\end{subfigure}
\begin{subfigure}{0.5\textwidth}
  \centering
  \includegraphics[width=.98\linewidth]{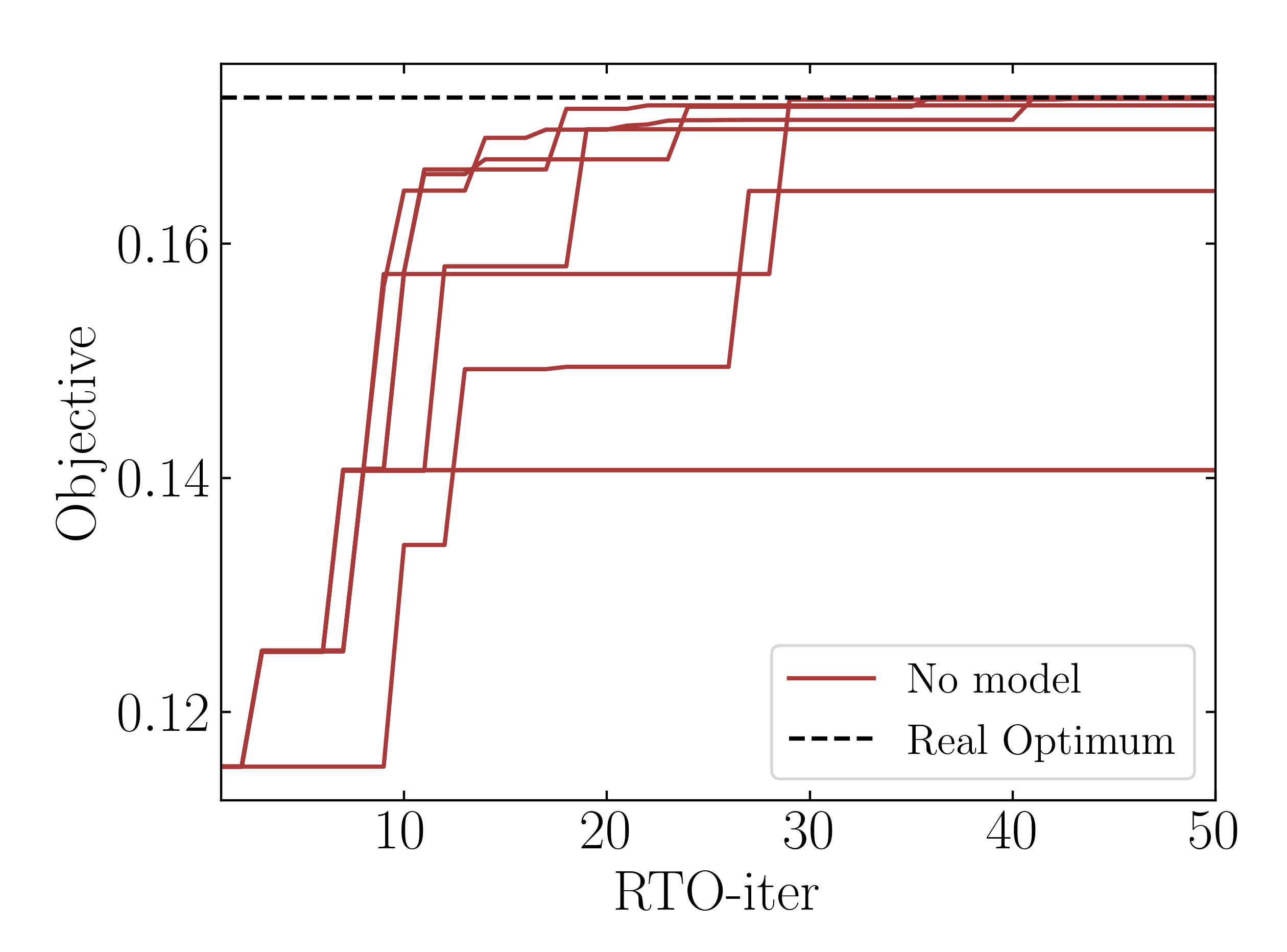}\vspace{-.5em}
  \caption{RTO iterates without prior model}
  \label{fig:PBR_noprior}
\end{subfigure}
\begin{subfigure}{0.5\textwidth}
  \centering
  \includegraphics[width=.98\linewidth]{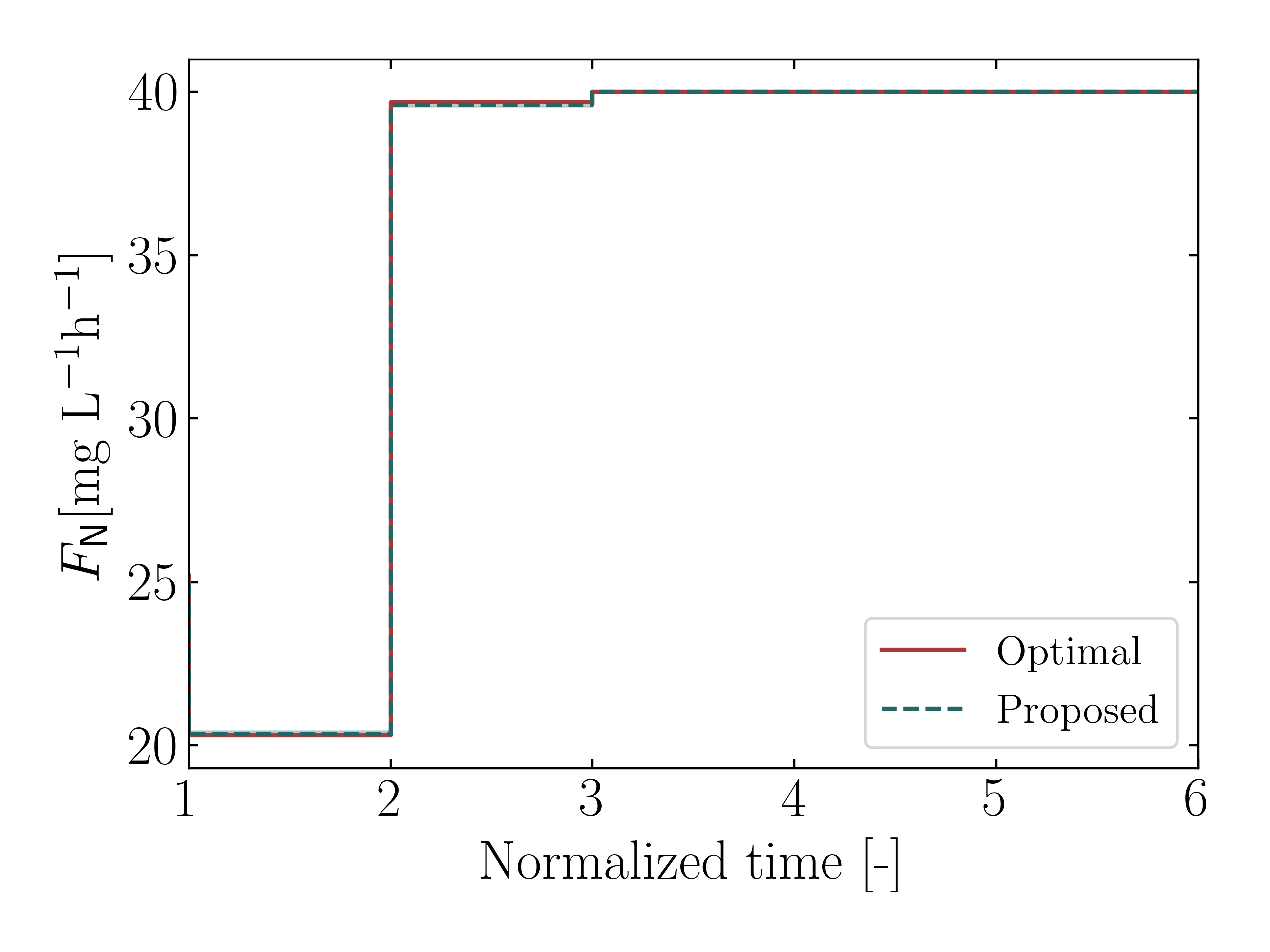}\vspace{-.5em}
  \caption{Optimized input $F_{\rm N}(t)$ at final iteration with prior model}
  \label{fig:PBR_FN}
\end{subfigure}
\begin{subfigure}{0.5\textwidth}
  \centering
  \includegraphics[width=.98\linewidth]{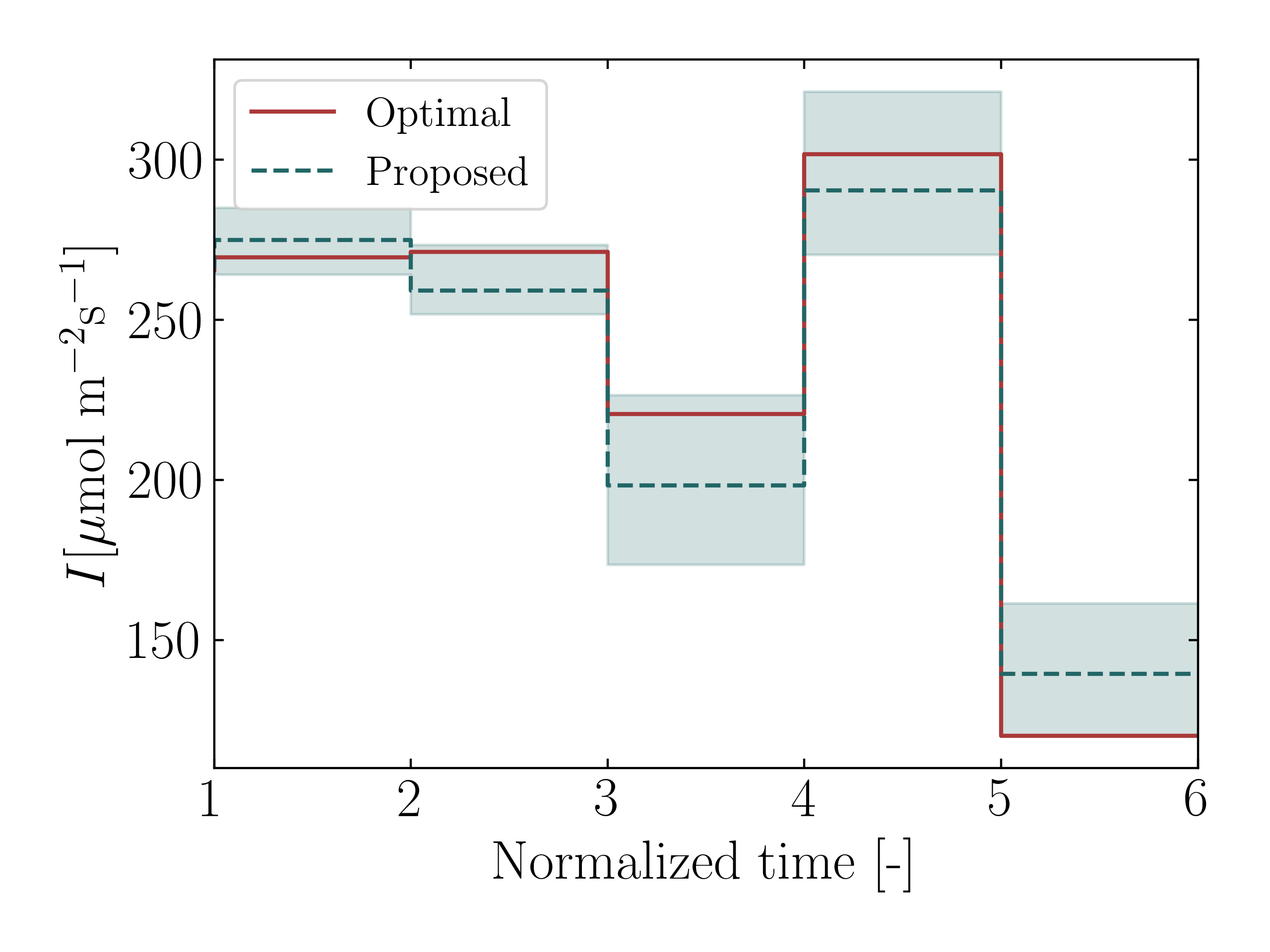}\vspace{-.5em}
  \caption{Optimized input $I(t)$ at final iteration with prior model}
  \label{fig:PBR_I}
\end{subfigure}
\caption{RTO results for the photobioreactor case study (Problem~\ref{eq:PBR}) using Algorithm~\ref{alg:MA-GP-TR}. {\bf (a)} Evolution of process cost with the RTO iterations for 8 process noise realizations with the nominal model \eqref{ode:Xmod}--\eqref{ode:Pmod} used as prior. {\bf (b)} Evolution of process cost with the RTO iterations for 8 process noise realizations without a prior model. Only the feasible iterates are shown. {\bf (c), (d)} Comparison between the optimal inputs $F_{\rm N}(t)$ and $I(t)$ and the RTO iterates after 50 iterations with the nominal model \eqref{ode:Xmod}--\eqref{ode:Pmod} used as prior; the envelopes are the same 8 noise realizations as in (a) and the dotted lines show one particular realization. \label{fig:PBR}}
\end{figure}

The initial GPs are trained with 13 feasible data points, which were obtained via trial-and-error, and the initial trust region encloses all of these points. All of the constraints are considered unrelaxable in Algorithm~\ref{alg:MA-GP-TR}. But unlike the other case studies, the trust-region radius is not reduced after back-tracking from an infeasible iterate as this was found to significantly hinder the progression of the RTO iterates.

The performance of Algorithm~\ref{alg:MA-GP-TR} with the nominal model \eqref{ode:Xmod}--\eqref{ode:Pmod} as prior and with the EI acquisition function is presented in Figure~\ref{fig:PBR_prior} for multiple realizations of the process noise. All of the runs are seen to reach a neighborhood of the plant optimum within 25--40 iterations, which may be considered fast given the large number of manipulated inputs. The optimized input profiles corresponding to $F_{\sf N}(t)$ and $I(t)$ after 50 iterations are shown in Figures~\ref{fig:PBR_FN} \& \ref{fig:PBR_I}, respectively, for the same noise realizations. It can be checked that all of these input profiles are indeed in excellent agreement with the plant optimum\deleted{, the larger variation range for the input $I(t)$ being attributed to its lower sensitivity compared to $F_{\sf N}(t)$}. \added{The smaller variation range for the input $F_{\sf N}(t)$ is attributed to the bang-bang nature of its optimal trajectory, which is thus determined by process constraints; while the optimal trajectory of the input $I(t)$ is comprised of interior arcs, which are known to be less sensitive \citep{Deshpande2012}.}

For comparison, the performance of the same algorithm without a prior model (model-free RTO) is reported in Figure~\ref{fig:PBR_noprior}. Notice that the behavior is now much more inconsistent across the various RTO runs, with certain runs converging to the plant optimum after just 20 iterations, while others failing to reach the plant optimum and remaining vastly suboptimal after 50 iterations. These results confirm that the use of a nominal model in the manner of a prior constitutes an effective derisking strategy in higher-dimensional RTO problems.

\section{Conclusions and Future Directions}
\label{sec:concl}

The main contribution of this paper lies in the development of an improved modifier-adaptation algorithm by integrating ideas from the related fields of Bayesian optimization and derivative-free optimization. On the one hand, trust-region techniques robustify the search by mitigating risk during the exploration \replaced{or accelerating the search whenever possible}{and enable building on an established convergence theory, e.g. for unconstrained RTO problems}. On the other hand, GPs are ideally suited to capture the plant-model mismatch or process noise in RTO, and a GP's variance estimator can drive the exploration by means of an acquisition function. \added{Special emphasis has been on algorithms that target good practical performance, rather than certifying global convergence.}

The \replaced{performance of the proposed algorithm has}{se benefits have} been analyzed and illustrated with numerical case studies, including a challenging batch-to-batch optimization problem with a dozen inputs and a large number of constraints. \added{Integrating an acquisition function in the modified optimization model provides clear benefits in terms of steering the iterates to the neighborhood of a plant optimum, especially in the presence of noise.} The paper has also investigated the benefits of embedding a prior (nominal) process model in the RTO scheme, instead of relying entirely on process data as in model-free RTO. The numerical case studies suggest that embedding a prior model can provide an effective derisking strategy against process noise. In practical applications, this added reliability could outweigh the benefits of model-free RTO, for instance in terms of ease of design and maintainability.

Future work will be geared towards improving the reliability of modifier-adaptation RTO schemes further, including the consideration of acquisition functions for the process constraints in order to promote exploration of the feasible region and accuracy of the GP surrogates\added{; and building on established convergence theory in the field of stochastic derivative-free optimization}. Another promising direction entails incorporating transient information to train the GPs, with a view to enabling dynamic real-time optimization. 


\section*{Acknowledgements}

This paper is based upon work supported by the UK Research and Innovation, and Engineering and Physical Sciences Research Council under grants EP/T000414/1 and EP/P016650/1. Financial support from Shell and FAPESP under grant 2014/50279-4, ANP and CNPq Brasil under grant 200470/2017-5) is gratefully acknowledged. This project has also received funding from the European Union's Horizon 2020 research and innovation programme under the Marie Sklodowska-Curie grant agreement No 675215.

\section*{Supporting Information}
The python codes implementing the numerical case studies can be retrieved from the following Git repository: \href{https://github.com/omega-icl/ma-gp}{https://github.com/omega-icl/ma-gp}. 



\bibliography{CACE_revision1}

\appendix
\gdef\thesection{\Alph{section}} 
\makeatletter
\renewcommand\@seccntformat[1]{Appendix \csname the#1\endcsname.\hspace{0.5em}}
\makeatother

\section{Global Convergence in Unconstrained RTO Problems}
\label{app:convergence}

This appendix \replaced{summarizes key global convergence results in derivative-free trust-region methods}{formalizes the global convergence properties of Algorithm~\ref{alg:MA-GP-TR} to a first-order critical point in unconstrained RTO problems}. \added{These convergence certificates are available for unconstrained optimization problems, though convergence in constrained optimization problems may also be certified, e.g., by converting them into unconstrained problems using penalty functions \citep{Larson2019}. In this context, the modifier-adaptation scheme in Algorithm~\ref{alg:MA-GP-TR} simplifies as follows: }\deleted{In this setup, the reduced gradient (Step~1) corresponds to the modified cost gradient $\boldsymbol{\nabla} [G_0+\mu^k_{\delta G_0}]$;} the acquisition function of the optimization subproblems (Step~2) is simply the modified cost $[G_0+\mu^k_{\delta G_0}]$; \deleted{the process measurements (Step 3) are noiseless;} and $n_g=0$, so \added{the reduced gradient (Step~1) corresponds to the modified cost gradient $\boldsymbol{\nabla} [G_0+\mu^k_{\delta G_0}]$ and} the feasibility test (Step~4) \replaced{becomes superfluous}{is not needed}. 

\replaced{In the idealized scenario of noiseless process measurements,}{Building on} established convergence theory from the field of derivative-free optimization \deleted{\citep[][Chapter~10]{Conn2009}, a set of sufficient conditions} relies on the following assumptions:

\begin{assumption}\label{ass:plant}
The process cost $G^{\rm p}_0$ is continuously differentiable with Lipschitz continuous gradient and bounded from below on the neighborhood $\bigcup_{{\bf u}\in\mathcal{U}} \mathcal{B}({\bf u};\Delta_{\rm max})$ of the input domain $\mathcal{U}$ for some radius $\Delta_{\rm max}>0$.
\end{assumption}

\begin{assumption}\label{ass:model}
The modified cost $[G_0+\mu^k_{\delta G_0}]$ is fully linear on $\mathcal{B}({\bf u}^k;\Delta^k)$ at every iteration $k=0, 1,\ldots$ of Algorithm~\ref{alg:MA-GP-TR}; that is, $[G_0+\mu^k_{\delta G_0}]$ is continuously differentiable with Lipschitz continuous gradient, and there exist global constants $\kappa_{\rm ef},\kappa_{\rm eg}<\infty$ (independent of $k$) such that:
\begin{align*}
\left\|\boldsymbol{\nabla} G^{\rm p}_0({\bf u}^k + {\bf d}) - \boldsymbol{\nabla}[G_0+\mu^k_{\delta G_0}] ({\bf u}^k + {\bf d})\right\| \ &\leq \kappa_{\rm eg}\ \Delta^k\\
\left|G^{\rm p}_0({\bf u}^k + {\bf d}) - [G_0+\mu^k_{\delta G_0}] ({\bf u}^k + {\bf d})\right| \ &\leq \kappa_{\rm ef}\ (\Delta^k)^2
\end{align*}
for all ${\bf d}\in \mathcal{B}({\bf u}^k;\Delta^k)$\deleted{ and all iterations $k$}.
\end{assumption}

\noindent
Under Assumptions~\ref{ass:plant} and \ref{ass:model}, \added{Theorem~10.13 in \citet{Conn2009} proves the convergence of} the \added{noisefree, unconstrained RTO} iterates produced by Algorithm~\ref{alg:MA-GP-TR} \deleted{for an unconstrained RTO problem with noiseless measurements are globally convergent} to a first-order critical point.
\replaced{The}{Besides standard smoothness and boundedness assumptions on the process and model cost functions, the} key assumption \added{here }is the need for a fully linear model (Assumption~\ref{ass:model})\added{, at least} on \deleted{every} iteration\added{s that do not yield a sufficient decrease in the objective value \citep{Conn2009b,Conn2009}}. This ensures that the \added{surrogate} model of the objective function has uniformly good local accuracy, similar in essence to the local behavior of first-order Taylor model. \added{Techniques for constructing fully-linear models are also well established. For instance, $\Lambda$-poised sets of points over a trust region can be generated for a wide variety of RBF-based surrogates \citep{Wild2008,Wild2013}, including Gaussian RBFs popularly used in GPs.}\deleted{In their multifidelity optimization algorithm, \citet{March2012} also enforced full linearity of the surrogate models on every iteration using a sampling policy developed by \citet{Wild2008}. Likewise, the constrained flowsheet optimization algorithm by \citet{Eason2016,Eason2018} posits full linearity of the surrogate models on every iteration.}

\added{In the more practical scenario of noisy process measurements, the surrogate models become probabilistic in nature. The global convergence of trust-region methods in such setup can be established based on the following extra or modified assumptions \citep{Bandeira2014}:}

\begin{assumption}\label{ass:noise}
\added{The additive process noise observed in measuring $G_0^{\rm p}$ is drawn from a distribution with mean zero and finite variance.}
\end{assumption}

\begin{assumption}\label{ass:randmodel}
\added{The modified cost $[G_0+\mu^k_{\delta G_0}]$ is fully linear with probability $\alpha$ on $\mathcal{B}({\bf u}^k;\Delta^k)$ for all sufficiently large iteration $k$; that is,}
\begin{align*}
\added{\mathbb{P}\left( \text{$[G_0+\mu^k_{\delta G_0}]$ is a fully-linear model on $\mathcal{B}({\bf u}^k;\Delta^k)$}~\middle|~\mathcal{F}_k \right) \geq \alpha}
\end{align*}
\added{where $\mathcal{F}_k$ is  the  filtration  of  the  random  process  up to  the  current iteration.}
\end{assumption}

\noindent
\added{Under Assumptions \ref{ass:plant}, \ref{ass:noise} and \ref{ass:randmodel} -- with additional conditions linking the probability level $\alpha$ to the parameters $\gamma_{\rm inc}$ and $\gamma_{\rm red}$ of the trust-region algorithm -- \citet[][Theorem~1]{Larson2016} proved the convergence in probability of a variant of Algorithm~\ref{alg:MA-GP-TR} to a first-order critical point. Under similar assumptions, \citet{Chen2018} proved almost sure convergence to a stationary point. The pivotal condition is again Assumption~\ref{ass:randmodel}, as a guarantee that the surrogate models have good accuracy with sufficiently high probability. A major benefit of the probabilistic full-linearity property is that it removes the need for $\Delta$-poised sets on every iteration. One way of satisfying this property is by regressing a sufficiently large number of sampling points, e.g. using linear regression \citep{Larson2016}. Recently, \citet{Shukla2020} proved that GP surrogates are probabilistically fully linear, thereby supporting their faster global convergence compared to local, linear or quadratic surrogate models.}

\deleted{Unfortunately, enforcing full linearity at every step of an RTO system would involve taking extra samples by perturbing the process in all directions around the current iterate, which is undesirable or even impractical. Though it should be noted that the global convergence property of derivative-free trust-region algorithms may still hold if the models are not fully linear on every iteration, provided that the iterates yield a sufficient decrease in the objective function \citep{Conn2009b,Conn2009}. Instead, the main idea in Algorithm~\ref{alg:MA-GP-TR} is to consider acquisition functions from the field of Bayesian optimization to promote exploration within the trust region (cf. Problem~\ref{eq:modified_problem_GP+TR}). This approach does not come with a formal convergence proof, but its efficiency and reliability can be established with numerical case studies.}


\linespread{1}

\begin{table*}[t]
\begin{framed}
\renewcommand{\nomname}{}
\printnomenclature
\end{framed}
\end{table*}
\pagestyle{empty}

\end{document}